\documentclass[12pt,twoside,reqno]{amsart}
\usepackage[colorlinks=true,citecolor=blue]{hyperref}
\usepackage{mathptmx, amsmath, amssymb, amsfonts, amsthm, mathptmx, enumerate, color,mathrsfs}
\usepackage{graphicx}

\usepackage{multirow}
\usepackage{epstopdf}
\usepackage{multicol}
\usepackage{algorithm}
\usepackage{algorithmic}
\usepackage{epstopdf}
\usepackage{amsmath,amsfonts,amssymb,mathrsfs,cases,mathdots,colordvi,url,extarrows,graphicx}
\usepackage{relsize}
\usepackage{mathrsfs}
\usepackage{color}
\usepackage{indentfirst}
\usepackage{float}
\usepackage[numbers,sort&compress]{natbib}

\newtheorem{defn}{Definition}[section]

\newtheorem{lemma}{Lemma}[section]
\newtheorem{thm}{Theorem}[section]
\newtheorem{prop}{Proposition}[section]
\newtheorem{cor}{Corollary}[section]
\newtheorem{example}{Example}[section]
\newtheorem{remark}{Remark}[section]

\renewcommand{\Box}{\rule{2.2mm}{2.2mm}}

\def\beginproof{\par\noindent {\bf Proof.}\ \ }
\def\endproof{\hskip .5cm $\Box$ \vskip .5cm}

\def\beginproof{\par\noindent {\bf Proof.}\ \ }
\def\endproof{\hskip .5cm $\Box$ \vskip .5cm}

\begin{document}
	\setcounter{page}{1}

	\vspace*{2.0cm}
	\title[subdifferential of the value function]
	{Directional subdifferential of the value function}
	\author[K. Bai, J. J. Ye]{ Kuang Bai$^{1}$, Jane J. Ye$^{2,*}$}
	\maketitle
	\vspace*{-0.6cm}

	\begin{center}
		{\footnotesize
			
			$^1$Department of Applied Mathematics, The Hong Kong Polytechnic University, Hong Kong, China\\
			$^2$Department of Mathematics and Statistics, University of Victoria, Canada\\
			Dedicated to Francis Clarke on the occasion of his seventy-fifth birthday
	}\end{center}

	\vskip 4mm {\footnotesize \noindent {\bf Abstract.}
		The directional subdifferential of the value function gives an estimate on how much the optimal value  changes under a perturbation in a certain direction.  In this paper we derive upper estimates for the  directional limiting and singular subdifferential of the value function for a very general parametric optimization problem. 
		We obtain a  characterization for the directional Lipschitzness of  a  locally lower semicontinuous function  in terms of the directional subdifferentials. Based on this characterization and the derived upper estimate for the directional singular subdifferential, we are able to  obtain a  sufficient condition for the directional Lipschitzness of the value function. Finally, we specify these results for various  cases  when all functions involved are smooth, when the perturbation is additive, when the constraint is independent of the parameter, or when the constraints are equalities and inequalities. Our results extend the corresponding  results on the sensitivity of the value function to allow directional perturbations. Even in the case of full perturbations, our results recover  or  even extend some existing results, including the Danskin's theorem.

		
		 }

	\renewcommand{\thefootnote}{}
	\footnotetext{ $^*$Corresponding author.
		\par
		E-mail address: kuang.bai@polyu.edu.hk (K. Bai), janeye@uvic.ca (J. J. Ye).
		\par
		Received xx, x, xxxx; Accepted xx, x, xxxx.

		\rightline {\tiny   \copyright  2022 Communications in Optimization Theory}}

\section{Introduction}
In this paper we consider a  parametric optimization problem in the form:
\begin{eqnarray*}
(P_x)~~~~\quad \min_{y\in\mathbb R^m} &&f(x,y)\\
\quad s.t.\ && P(x,y)\in \Gamma,
\end{eqnarray*}
where $x\in\mathbb R^n$ denotes the parameter, $f:\mathbb R^{n+m}\rightarrow\mathbb R,\ P:\mathbb R^{n+m}\rightarrow\mathbb R^p$ and $\Gamma\subseteq\mathbb R^p$ is a closed set. 
Unless otherwise stated, we assume that all functions $f, P$ are locally Lipschitz continuous.
 Such a problem is very general. In the non-parametric case, it is sometimes referred to as the mathematical program with geometric constraints; see e.g. \cite{GYZ} or as a set-constrained optimization problem; see e.g.  \cite{GYZnew}. In the case where $\Gamma=\mathbb{R}^{p_1}_-\times \{0\}^{p_2}$ with $p_1+p_2=p$, the parametric optimization problem $(P_x)$ becomes  a parametric nonlinear programming problem with equality and inequality constraints and in the case where $\Gamma$ is a convex cone, it is a parametric conic program as studied in \cite{BS}.

In practice, it is important to know how  the optimal value of an optimization problem  changes subject to perturbation.
For this purpose, it is interesting to study certain properties such as Lipschitz continuity and differentiability of the associated (optimal) value function/marginal function
$$V(x):=\inf_y \{f(x,y)| P(x,y)\in \Gamma\},$$ 
where by convention, $V(x)$ is defined to be $+\infty$ if the feasible map
$${\mathcal F}(x):=\{y|P(x,y)\in \Gamma \}$$ is empty at $x$. We define the solution map of $(P_x)$  by
$$S(x):=\arg\min_y\{ f(x,y)| P(x,y)\in \Gamma\}.$$
Sensitivity analysis of the value function consists of the study of its directional differentiability and  subdifferentials. In this paper we mainly concern about subdifferentials of the value function and refer the reader to the topic on the directional differentiability of the value function in a forthcoming paper \cite{BY2022a}.

Since it is known that in general the value function is not smooth even when all problem data are smooth, in the literature, one usually tries to  give upper estimates of certain  subdifferentials of the value function and then use them to obtain some useful information.
The classical results in this regard were given for smooth nonlinear programs with equality and inequality constraints.   Gauvin and Dubeau \cite[Theorems 5.1 and 5.3]{GD} showed that for a smooth nonlinear program, if 
$S(\bar x)\not =\emptyset$, the solution map $S(x)$ is uniformly compact near $\bar x$, and Mangasarian-Fromovitz constraint qualification (MFCQ) holds at each $\bar y\in S(\bar x)$ then $V(x)$ is Lipschitz continuous near $\bar x$ and the following upper estimate holds for the Clarke subdifferential/generalized gradient of the value function:
\begin{equation}\label{eqn1}
\partial^c V(\bar x) \subseteq {\rm co} \bigcup_{\bar y\in S(\bar x) } \left \{ 
 \nabla_x L(\bar x,\bar y,\lambda)\left | \lambda \in \Lambda(\bar x,\bar y)  \right . \right  \},\end{equation}
where $L(x,y,\lambda):=f(x,y)+P(x,y)^T\lambda$ is the Lagrange function and $\Lambda(x,y)$ is the set of the Lagrange multipliers for problem $(P_{ x})$ at $y\in S(x)$. Moreover by \cite[Corollary 5.4]{GD}, if MFCQ is replaced by the linear independence constraint qualification  in the above, then $V(x)$ is Clarke regular and (\ref{eqn1}) holds as an  equality. Furthermore if the solution is unique, then the value function is  smooth.  Clarke \cite{Clarke} considered an additively (right-hand side) perturbed  nonsmooth optimization problem with equality, inequality and an abstract constraint.
The restricted inf-compactness condition (see \cite[Hypothesis 6.5.1]{Clarke}, \cite{GLYZ}) which is one of the weakest sufficient conditions for the lower semicontinuity of the value function was introduced and  upper estimates not only for the generalized gradient but also for the asymptotic generalized gradient were  obtained in Clarke \cite{Clarke}. Lucet and Ye \cite{LY,LYErru} gave
upper estimates for the limiting subdifferential and the singular subdifferential of the value function for a perturbed  nonsmooth optimization problem with equality, inequality, an abstract constraint and a variational inequality constraint.
For nonsmooth nonlinear programs with  equality and inequality constraints, Ye and Zhang \cite{YZ}   obtained upper estimates of the limiting and the singular subdifferential in terms of the enhanced multipliers and the abnormal enhanced multipliers respectively. Since the set of the enhanced multipliers is contained in the set of  the standard multipliers, the obtained upper estimates in \cite{YZ} are  sharper than those in terms of the standard multipliers. The results in \cite{YZ} have been extended by  Guo et al. \cite{GYZ} to the parametric program with geometric constraints, and  to the  parametric mathematical program with equilibrium constraints (MPEC) by  Guo et al. in \cite{GLYZ}.

Sometimes, there are advantages in considering  perturbations only in a certain direction. 
To deal with these kinds of requirements, a directional version of the  limiting normal cone and subdifferential  have been introduced independently by Ginchev and Mordukhovich \cite{GM} and Gfrerer \cite{Gfr13}. The directional limiting normal cone and subdifferential are in general   
 smaller than their non-directional counterparts and possess rich calculus 
  (see {\cite{BHA,BP2022,Long,TC2017}}). These directional objects allow one to consider directional optimality conditions which are sharper than the nondirectional one (see e.g Gfrerer \cite{Gfr13} and Bai and Ye \cite{BY}), directional constraint qualifications which are weaker than its nondirectional counterparts (see e.g., \cite{Gfr13,BYZ,Gfr13b,BP2022A,BP2022E,BP2022F,Hcor}), optimality conditions for nonconvex mathematical programs (see e.g., {\cite{BHYZ,Gfr14,GYZnew,TC2017}}) and stability analysis of constraint systems/set-valued maps (see e.g., {\cite{B2021,GO2016,GO2016b,YZ17}}).
 Recently Bai and Ye \cite{BY} introduced a directional version of the Clarke subdifferential and gave upper estimates for the directional limiting/Clarke subdifferential of the value function for  parametric smooth nonlinear programs with equality and inequality constraints.
 
 The upper estimate for the  directional limiting subdifferential of the value function  in  Bai and Ye \cite[Theorem 4.2]{BY} was obtained under several assumptions including the relaxed constant rank regularity (RCR-regularity) condition and the Robinson stability (RS), and the results are only for smooth nonlinear programs.  In this paper we obtain upper estimates not only for the directional limiting subdifferential but also for the singular directional limiting subdifferential of the value function. Our results are obtained for the general parametric optimization problem $(P_x)$ under only the metric subregularity/calmness condition. 
 For the  additive perturbation, our results extend and recover the classical results in Clarke \cite{Clarke}.  In the case of the smooth nonlinear program, our results do not need the RCR-regularity and the RS conditions as in   Bai and Ye \cite[Theorem 4.2]{BY}. Moreover we extend  the celebrated Danskin's theorem to a directional version.

 We organize the paper as follows. In the next section, we provide the notation and preliminary results. In section 3, we derive upper estimates for the directional subdifferentials  of the value function and give new sufficient conditions for the directional Lipschitz continuity of the value function.  Finally, in section 4, we apply the results to various special cases.

\section{Preliminaries and preliminary results} We first give notation that will be used in the paper.  Let $\Omega$ be a set. By $x^k\xrightarrow{\Omega}\bar{x}$ we mean $x^k\rightarrow\bar{x}$ and for each $k$, $x^k\in \Omega$.  By $x^k\xrightarrow{u}\bar x$, we mean that the sequence $\{x^k\}$ approaches $\bar x$ in direction $u$, i.e., there exist $t_k\downarrow 0, u^k\rightarrow u$ such that $x^k=\bar x+ t_k u^k$. By $f(t)=o(t)$, we mean that $f(t)$ is a function such that $\lim_{t\downarrow0}\frac{f(t)}{t}=0$. $\mathbb S$ denotes the unit sphere.
$\mathbb B$ denotes the unit open ball and $\mathbb B_\sigma (\bar x)$ denotes the open ball centered at $\bar x$ with radius equal to $\sigma$.
${\mathcal B}$ denotes the  closed unit ball. For any  $x\in \mathbb R^n$, we denote by $ \|x\|$ the Euclidean norm.
 For a set $\Omega$, we denote by  ${\rm co}\Omega$, ${\rm clco}\Omega$, $\overline{\Omega}$,  bd$\Omega$ and $\Omega^\perp$  its convex hull, its closure of the convex hull, its closure, its boundary and its orthogonal complement, respectively.  By  ${\rm dist}(x,\Omega):=\inf\{\|x-y\||y\in\Omega\}$, we denote the distance from a point $x$ to set $\Omega$. For a single-valued map $\phi:\mathbb R^n\rightarrow\mathbb R^m$, we denote by $\nabla \phi(x)\in \mathbb{R}^{m\times n}$   the Jacobian matrix of $\phi$  at $x$ and for a function $\varphi:\mathbb R^n\rightarrow\mathbb R$, we denote by $\nabla \varphi(x)$ both the gradient and the Jacobian of $\varphi$  at $x$. We denote the extended real line by  $\overline{\mathbb{R}}:=[-\infty,\infty]$ and for an extended-valued function $f:\mathbb{R}^n\rightarrow \overline{\mathbb{R}}$, we define the effective domain of $f$ as $ \mbox{dom}f:=\{x|f(x) <\infty\}$. For a set-valued map  $\Phi:\mathbb R^n\rightrightarrows\mathbb R^m$, we define its graph by ${\rm gph}\Phi:=\{(x,y)| y\in \Phi(x)\}$ and its domain by ${\rm dom}\Phi:=\{x|  \Phi(x)\not =\emptyset\}$.

When $d=0$ the following  definition coincides with the Painlev\'e-Kuratowski inner/lower and outer/upper  limit of $\Phi$ as $x\rightarrow\bar x$ respectively.
\begin{defn} Given a set-valued map $\Phi:\mathbb R^n\rightrightarrows\mathbb R^m$ and a direction $
d \in\mathbb R^n$, the inner/lower and outer/upper limit of $\Phi$ as $x\xrightarrow{d}\bar x$ respectively is defined by
	\begin{align*}
	\liminf_{x\xrightarrow{d}\bar x} \Phi(x):=\{y\in\mathbb R^m|&\forall \ \mbox{sequences}\ t_k\downarrow0, d^k\rightarrow d, \exists y^k\rightarrow y \mbox{ s.t. } y^k\in \Phi(\bar x+t_kd^k)\},\\
	\limsup_{x\xrightarrow{d}\bar x} \Phi(x):=\{y\in\mathbb R^m|&\exists   \ \mbox{sequences}\ t_k\downarrow0, d^k\rightarrow d,  y^k\rightarrow y \mbox{ s.t. } y^k\in \Phi(\bar x+t_kd^k)\},
	\end{align*}
	respectively.
\end{defn}
We review the various concepts of tangent and normal cones below (see, e.g., \cite{Clarke,Danskin}, \cite[Definitions 6.1 and 6.3]{RW}, \cite[Theorem 3.57]{Aub2} and \cite[Definition 2.54]{BS}).
\begin{defn}[Tangent Cones and Normal Cones] 
	Given a set $\Omega\subseteq\mathbb{R}^n$ and a point $\bar{x}\in \Omega$, the inner tangent cone to $\Omega$ at $\bar{x}$ is defined as
	$$T^i_\Omega(\bar{x}):=
	\left \{d\in\mathbb{R}^n|\forall t_k\downarrow0, \exists d_k\rightarrow d\ \mbox{ s.t. } \bar{x}+t_kd_k\in \Omega\ \forall k\right \},$$
	the tangent/contingent  cone to $\Omega$ at $\bar{x}$ is defined as
	$$T_\Omega(\bar{x}):=
	\left \{d\in\mathbb{R}^n|\exists t_k\downarrow0, d_k\rightarrow d\ \mbox{ s.t. } \bar{x}+t_kd_k\in \Omega\ \forall k\right \}.$$	
	The regular normal cone,  the limiting normal cone and the Clarke normal cone to $\Omega$ at $\bar{x}$ can be defined as 
\begin{eqnarray*}\widehat{N}_\Omega(\bar{x})&:=&\left \{ \zeta\in \mathbb{R}^n\bigg| \langle \zeta ,x-\bar{x}\rangle \leq o(\|x-\bar x\|) \quad \forall x\in \Omega \right \},\\
	 N_\Omega(\bar{x})&:=&
	 \left \{\zeta\in \mathbb{R}^n\bigg| \exists \ x_k\xrightarrow{\Omega}\bar{x},\ \zeta_k{\rightarrow}\zeta\ \text{such that}\ \zeta_k\in\widehat{N}_\Omega(x_k)\ \forall k\right 
	 \},\\
	N^c_\Omega(\bar{x})&:=& {\rm clco} N_\Omega(\bar{x}),  \end{eqnarray*}
	respectively.
	\end{defn}	
We say that a set $\Omega$ is geometrically derivable at $\bar x\in\Omega$ if $T^i_\Omega(\bar x)=T_\Omega(\bar{x})$. It is known that not only convex sets, but also polyhedral sets  (union of finitely many convex polyhedral sets) as well as a second-order complementarity set \cite[Proposition 5.2]{YZ17} 
are geometrically derivable. 
	
	
%
	\begin{defn}[Directional Normal Cone] \cite{GM,Gfr13}
    Given a set $\Omega \subseteq \mathbb{R}^n$, a point $\bar x \in \Omega$ and a direction $d\in\mathbb{R}^n$, the limiting normal cone to $\Omega$ at $\bar{x}$ in direction $d$ is defined by
    $$N_\Omega(\bar{x};d):=\left \{\zeta \in \mathbb{R}^n\bigg| \exists \ t_k\downarrow0, d_k\rightarrow d, \zeta_k\rightarrow\zeta  \mbox{ s.t. } \zeta_k\in \widehat{N}_\Omega(\bar{x}+t_kd_k)\ \forall k \right \}.$$
\end{defn}
It is obvious that $N_{\Omega}(\bar x; 0)=N_{\Omega}(\bar x)$, $N_{\Omega}(\bar x; d)=\emptyset$ if $d \not \in T_\Omega(\bar x)$ and $N_{\Omega}(\bar x;d)\subseteq N_\Omega(\bar x)$. It is also obvious that for all $d\in T_\Omega(\bar x) \setminus T_{{\rm bd} \Omega}(\bar x)$, one has $N_\Omega(\bar x;d)=\{0\}$.
Moreover when $\Omega$ is convex, by \cite[Lemma 2.1]{Gfr14} the directional and the classical normal cone have the following relationship
\begin{equation}
N_\Omega(\bar x;d)=N_\Omega(\bar x)\cap \{d\}^\perp  \qquad \forall d\in T_\Omega(\bar x).\label{convNormal}
\end{equation}

In fact the regular normal cone in the definition of the directional normal cone can be equivalently replaced by the limiting normal cone as in the following proposition. In another word, the directional limiting normal cone has the outer semicontinuity property.
\begin{prop}\cite[Proposition 2]{GYZnew}\label{noc} Given a set $\Omega \subseteq \mathbb{R}^n$, a point $\bar x \in \Omega$ and a direction $d\in\mathbb{R}^n$, one has
$$N_\Omega(\bar{x};d)=\left \{ \zeta \in \mathbb{R}^n\big| \exists \ t_k\downarrow0, d_k\rightarrow d, \zeta_k\rightarrow\zeta  \mbox{ s.t. } \zeta_k\in {N}_\Omega(\bar{x}+t_kd_k)\ \forall k  \right \}.$$
\end{prop}

{Similarly, many classical concepts in variational analysis can have their directional versions. {To this end}, the following concept of a directional neighborhood is needed.}
\begin{defn}[Directional Neighborhood]\label{dn} (\cite[formula (7)]{Gfr13}).
  Given a direction $d\in \mathbb{R}^n$, and positive numbers $\varepsilon,\delta>0$, the directional neighborhood of direction $d$ is a set defined by 
\begin{equation*}
{\mathcal V}_{\varepsilon,\delta}(d):=\{z\in\varepsilon\mathbb{B}|\big\| \| d\| z-\| z\| d\big\|\leq \delta \|z\| \|d\|\}.
\end{equation*}
\end{defn}
By definition, it is clear that
$$
{\mathcal V}_{\varepsilon,\delta}(d):=\left \{ \begin{array}{ll}
\{0\} \cup \left \{ z \in \varepsilon\mathbb{B}\setminus \{0\} \mid \|\frac{z}{\|z\|} -\frac{d}{\|d\|} \|\leq \delta \right \} & \mbox{ if } d\not =0,\\
\varepsilon\mathbb{B} &  \mbox{ if } d =0 \end{array} \right . .
$$
Hence if the direction $d=0$, then the directional neighborhood is nothing but an open ball $\varepsilon\mathbb{B}$ while if the direction $d\not =0$, the directional neighborhood is a section of the open unit ball with the central angle determined by $\delta$.  The colored section in Figure 1 shows the graph of $\mathcal V_{\varepsilon,\delta}(d)$ when $d\neq0$. Note that unless $d=0$, a directional neighborhood is not an open set.
\begin{figure}[h!]
	\centering
	\includegraphics[height=5cm]{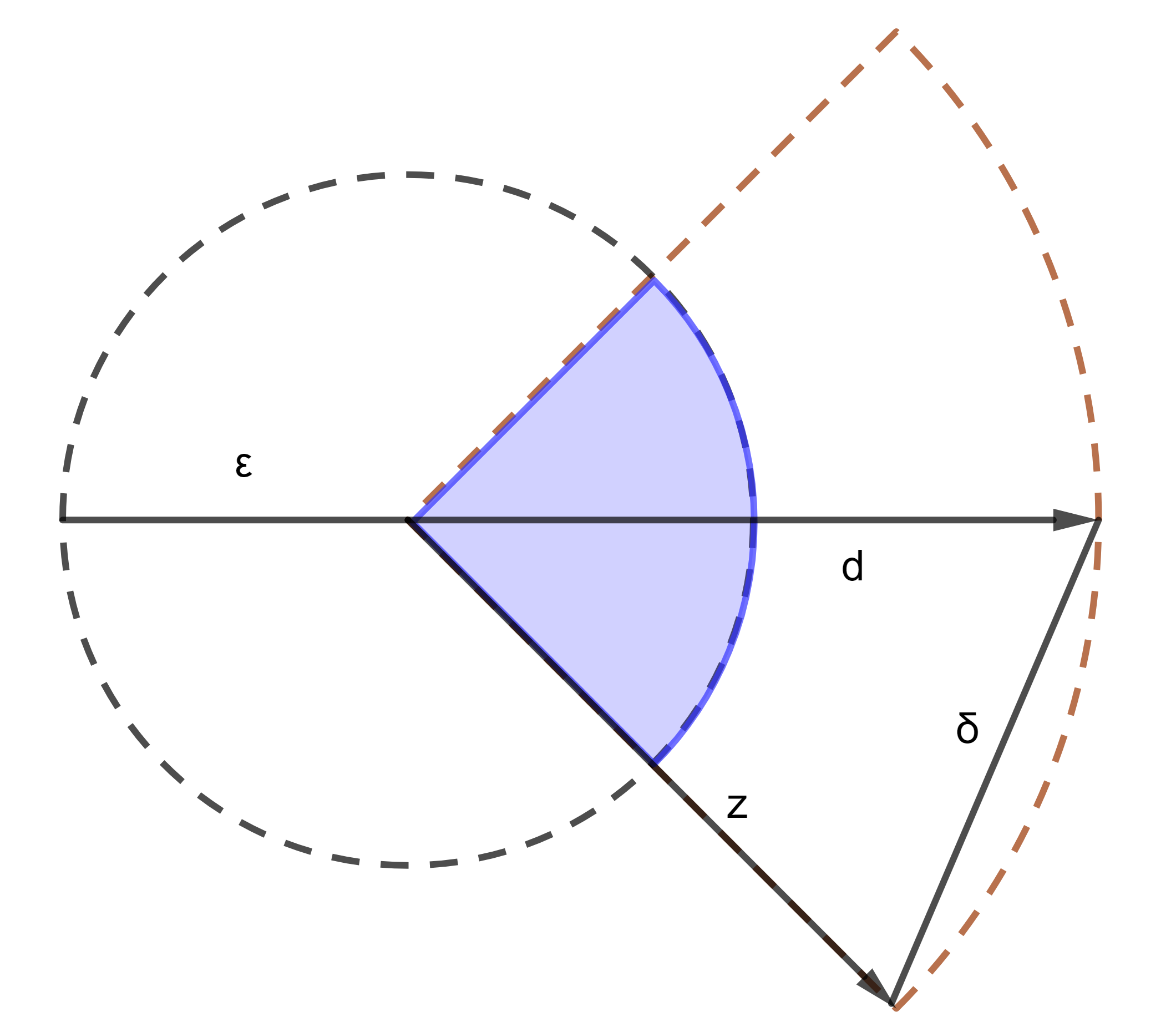}
	\caption{Directional neighborhood}
\end{figure}


With the directional neighborhood at hand, many classical concepts have been exteneded to their directional versions, and a directional version of the variational analysis has been stimulated (see discussions in the Introduction). 
 In the sequel, we list some useful results from directional variational analysis.
First we define a concept of the directional lower semicontinuity and continuity. Obviously these concepts can also be defined by replacing the neighborhood by the directional neighborhood in the classical concepts.
\begin{defn}[Directional Lower Semicontinuity and Continuity]\label{dlsc} Let $f:\mathbb{R}^n \rightarrow \overline{\mathbb{R}}$ be finite at $\bar x$. We say $f$ is lower semicontinuous (l.s.c.) at $\bar x$ in direction $u$ if 
$$f(\bar x) \leq \liminf_{x\xrightarrow{u}\bar x}f(x). $$ 
We say $f$ is continuous at $\bar x$ in direction $u$ if 
$$f(\bar x) = \lim_{x\xrightarrow{u}\bar x}f(x). $$ 
\end{defn}

Using the concept of the directional neighborhood, \cite{BHA} defined a directional version of the Lipschitz continuity for functions.
\begin{defn}[Directional Lipschitz Continuity](\cite[Page 719]{BHA})  We say that a single-valued mapping $\varphi (x):\mathbb R^n\rightarrow\mathbb R^m$ is  Lipschitz continuous around $\bar x$ in direction $u$ if there exists a scalar $L\geq 0$ and a directional neighborhood ${\mathcal V}_{\varepsilon, \delta}(u)$ of $u$ such that
		\[ \|\varphi(x)-\varphi(x')\|\leq L\|x-x'\| \quad \forall x,x'\in \bar x+{\mathcal V}_{\varepsilon, \delta}(u).
		\]
If $u=0$ in the above definition, we say $\varphi(x)$ is  Lipschitz continuous  around $\bar x$.
\end{defn}

\begin{defn}[Directional Derivatives]\label{derivatives}
Let $\phi:\mathbb R^n\rightarrow\mathbb{R} $.
The Dini upper directional derivative of $\phi$ at $x$ in the direction $u$ is 
$$\phi'_+(x;u):=\limsup_{t\downarrow0}\frac{\phi(x+tu)-\phi(x)}{t}.$$
The Dini lower directional derivative of $\phi$ at $x$ in direction $u$ is
$$
\phi'_-(x;u):=\liminf_{t\downarrow0}\frac{\phi(x+tu)-\phi(x)}{t}.
$$
The usual directional derivative of $\phi$ at $x$ in the direction $u$ is 
$$\phi'(x;u):=\lim_{t\downarrow0}\frac{\phi(x+tu)-\phi(x)}{t}$$
when this limit exists.  The directional derivative of $\phi$ at $x$ in direction $u$ in Hadamard sense  is
$$
	\phi'_H(x;u):=\lim_{k\rightarrow\infty}\frac{\phi(x+t_ku^k)-\phi(x)}{t_k}
$$
when this limit exists for all sequences $t_k\downarrow0$, $u^k\rightarrow u$. 
\end{defn}

\begin{defn}[Graphical derivative]
	Let $\phi:\mathbb R^n\rightarrow\mathbb{R}^m $.
	The graphical derivative of  $\phi$ at $x$ in direction $u$ is
	\[
	D\phi(x)(u):=\left\{
	d\left|\exists t_k\downarrow0, u^k\rightarrow u\ \mbox{s.t.}\ \lim_{k}\frac{\phi(x+t_ku^k)-\phi(x)}{t_k}=d
	\right.\right\}.
	\]
\end{defn}
When $\phi$ is directionally differentiable at $\bar x$ in direction $u$ in Hadamard sense, {$\phi(x)$ is also continuous at $\bar x$ in direction $u$, and the graphical derivative $D\phi(\bar x)(u)$, is equal to the singleton $\{\phi_H'(\bar x;u)\}$.}
It is easy to see that if $\phi$ is Lipschitz continuous around $x$ in direction $u$, then $\phi'(x;u)=\phi'_H(x;u)$. Furthermore, if $\phi(x)$ is also directionally differentiable, then 
$$ D\phi(x)(u)=\phi'(x;u)=\phi'_H(x;u)=\phi'_-(x;u)=\phi'_+(x;u).$$

We review the definition of some subdifferentials below.
\begin{defn}[Subdifferentials](\cite[Definition 8.3]{RW} and \cite{Clarke})\label{sd} 
Consider an extended-real-valued function $f:\mathbb{R}^n \rightarrow \overline{\mathbb{R}}$ and a point $\bar x\in  {\rm dom} f.$
The Fr\'{e}chet (regular) subdifferential of f at $\bar{x}$ is the set
\begin{eqnarray*}
\widehat\partial f(\bar{x}):=\left\{\xi\in\mathbb{R}^n|f(x)\geq f(\bar x)+\langle\xi,x-\bar x\rangle+o(\|x-\bar x\|)\right\},
\end{eqnarray*}
the limiting/Mordukhovich/basic  subdifferential of f at $\bar{x}$ is the set
\begin{align*}
\partial f(\bar{x}):=
\{\xi\in\mathbb{R}^n|\exists x^k\rightarrow  \bar x,  \ \xi^k\rightarrow\xi \ \mbox{ s.t. }  f(x^k)\rightarrow f(\bar x),\ \xi^k\in \widehat\partial f(x^k)\},
\end{align*}
the singular/horizon subdifferential of f at $\bar{x}$ is the set
\begin{align*}
\partial^\infty f(\bar{x}):=
\{\xi\in\mathbb{R}^n|\exists x^k\rightarrow  \bar x,  \tau_k\downarrow0 \mbox{ s.t. } \tau_k\xi^k\rightarrow\xi ,\   f(x^k)\rightarrow f(\bar x),\ \xi^k\in \widehat\partial f(x^k)\}.
\end{align*}
If $f(x)$ is  Lipschitz continuous around $\bar x$, then the Clarke subdifferential of $f$ at $\bar x$ can be equivalently defined as
\[
\partial^cf(\bar x):=co(\partial f(\bar x)).
\]
\end{defn}

The following concepts are the directional version of the limiting and the singular subdifferential studied in \cite{Gfr13,GM,Long}. Recently \cite{BHA} introduced a concept of the  limiting and the singular subdifferential in  directions not only from $\bar x$ but also from   $f(\bar x)$.
\begin{defn}[Directional Subdifferentials]\label{ads}  Let $f:\mathbb{R}^n\rightarrow \overline{\mathbb{R}}$ and $\bar x\in {\rm dom}f$. 
The  limiting subdifferential of $f$ at $\bar x$ in direction 
$u\in \mathbb{R}^n$ is defined as 
	\begin{eqnarray*}
	\partial f(\bar x;u):=\{\xi\in\mathbb{R}^n|\exists x^k\xrightarrow{u}\bar{x}, \xi^k\rightarrow\xi \mbox{ s.t. } f(x^k)\rightarrow f(\bar x), \   \xi^k\in \widehat{\partial} f( x^k)\},
	\end{eqnarray*} 
the  singular subdifferential of $f$ at $\bar x$ in direction 
$u\in \mathbb{R}^n$ is defined as 
\begin{eqnarray*}
	\partial^\infty f(\bar x;u):=\left \{\xi\in\mathbb{R}^n|\begin{array}{l}
	\exists x^k\xrightarrow{u}\bar{x}, \tau_k\downarrow 0, \tau_k\xi^k\rightarrow\xi
	 \mbox{ s.t. } f(x^k)\rightarrow f(\bar x), \   \xi^k\in \widehat{\partial} f(x^k)\end{array} \right \}.
\end{eqnarray*} 
\end{defn}

Furthermore, if $f(x)$ is Lipschitz continuous around $\bar x$ in direction $u$, \cite{BY} defined the directional Clarke subdifferential of $f$ at $\bar x$ in direction $u$ as
\[
\partial^cf(\bar x;u):=co(\partial f(\bar x;u)).
\]
It is clear that $\partial^c f(\bar x;0)=\partial^c f(\bar x)$. Moreover, \cite[Proposition 2.1]{BY} proved that \begin{equation}
\partial^cf(\bar x;u)=co\limsup_{x\xrightarrow{u}\bar x}\partial^cf(x),\ \partial^c(-f)(\bar x;u)=-\partial^cf(\bar x;u). \label{clarke} \end{equation}

Similar to   the directional limiting normal cone in Proposition \ref{noc}, in the definition of the  directional limiting and singular subdifferentials, one can replace the regular subdifferential by the limiting subdifferential.
The result for the directional limiting subdifferential is given in \cite[Theorem 5.4]{Long}. In the following proposition, we prove the result for the directional singular subdifferential.
\begin{prop}\label{adsoc}
	Let   $u\in\mathbb R^n$ and $f:\mathbb R^n\rightarrow\overline{\mathbb{R}}$ 
	with   $\bar x\in {\rm dom}f$.
	Then
	\begin{equation*}
	\partial f(\bar x;u)=\{\xi\in\mathbb{R}^n|\exists x^k\xrightarrow{u}\bar{x}, \xi^k\rightarrow\xi \mbox{ s.t. } f(x^k)\rightarrow f(\bar x), \   \xi^k\in {\partial} f(x^k)\},
	\end{equation*}
and	\begin{equation}\label{s}
	\partial^\infty f(\bar x;u)=\left \{\xi\in\mathbb{R}^n|\begin{array}{l}
	\exists x^k\xrightarrow{u}\bar{x},\tau_k\downarrow0,  \tau_k\xi^k\rightarrow\xi 
	\mbox{ s.t. } f(x^k)\rightarrow f(\bar x), \   \xi^k\in {\partial} f(x^k)\end{array} \right \}.
	\end{equation}
\end{prop}
\beginproof 
We  prove the one for the directional singular subdifferential.
Denote the right hand side of (\ref{s}) by RHS.
By definition of the directional subdifferential, it is easy to see that  $\partial^\infty  f(\bar x;u)\subseteq {\rm RHS}$. To show the reverse inclusion, consider any $\xi\in{\rm RHS}$. Then there exist sequences $t_k\downarrow0, \tau_k\downarrow 0, u^k\rightarrow u$ and $\xi^k\in\partial f(\bar x+t_ku^k)$ with $\tau_k\xi^k\rightarrow\xi,\ f(\bar x+t_ku^k)\rightarrow f(\bar x)$. For each $k$, by Definition \ref{sd}, there exists $(\tilde x^k,f(\tilde x^k),\tilde\xi^k)\in (\bar x+t_ku^k,f(\bar x+t_ku^k),\xi^k)+t^2_k\mathbb B^{n+1+n}$ such that $\tilde{\xi}^k\in\widehat{\partial}f(\tilde x^k)$. Then we have \[
\lim_{k\rightarrow \infty}\frac{\tilde x^k-\bar x}{t_k}=\lim_{k\rightarrow \infty}\frac{\tilde x^k-(\bar x+t_ku^k)}{t_k}+\lim_{k\rightarrow \infty}\frac{\bar x+t_ku^k-\bar x}{t_k}=u, 
\]
and 
$$\lim_{k\rightarrow \infty}f(\tilde x^k)= \lim_k(f(\bar x+t_ku^k)+o(t_k))=f(\bar x),\ \lim_{k\rightarrow \infty}\tau_k \tilde\xi^k=\lim_{k\rightarrow \infty}\tau_k\left(\xi^k+o(t_k)\right)=\xi.$$ It follows that  $\xi\in\partial^\infty f(\bar x;u)$. This shows that $\partial^\infty  f(\bar x;u)\supseteq {\rm RHS}$ and the proof is complete.
\endproof

Recall that an extended-real-valued function $f$ is said to be locally l.s.c. at $\bar x$ where $f(\bar x)$ is finite if and only if ${\rm epi} f$ is locally closed at $(\bar x, f(\bar x))$ \cite[Excercise 1.34]{RW}, which means that there exists $\varepsilon>0$ such that 
${\rm epi} f \cap \overline{\mathbb{B}}_\varepsilon(\bar x, f(\bar x))$ is closed. We now extend this   definition to a directional local lower semicontinuity.
\begin{defn}[Directional Local Lower Semicontinuity]\label{locallsc} Let $f:\mathbb R^n\rightarrow\overline{\mathbb R}$ with $\bar x\in {\rm dom}f$. Then $f(x)$ is said to be locally l.s.c.~at $\bar x$ in direction $u$ if  there exist positive scalars $\varepsilon,\delta$ such that epi$f\cap \overline{\bar x+{\mathcal V}_{\varepsilon,\delta}(u)}\times \overline{\mathbb{B}}_\varepsilon(f(\bar x))$ is closed.
\end{defn}
{The concept of the  directional local lower semicontinuity above is stronger than the directional  lower semicontinuity defined in Definition \ref{dlsc}.
We now give an example to show that the directional continuity does not imply the directional local lower semicontinuity.
\begin{example}
	consider the function
	\begin{align*}
		f(x)=\left\{\begin{array}{rl}
			x,\ & \mbox{if}\ x\ \mbox{is a positive rational number},\\
			-x,\ & \mbox{if}\ x\ \mbox{is a negative rational number},\\
			0,\ & \mbox{otherwise}.
		\end{array}\right.
	\end{align*}
Let $\bar x=0$. Since $\displaystyle \lim_{x\xrightarrow{1}\bar x}f(x)=f(\bar x)$, $f(x)$ is directionally continuous at $\bar x$ along direction $u=1$. However, epi$f\cap \overline{\bar x+{\mathcal V}_{\varepsilon,\delta}(u)}\times \overline{\mathbb{B}}_\varepsilon(f(\bar x))=$ epi$f\cap ([0, \varepsilon]\times [-\varepsilon,\varepsilon])$ is not closed for any $\varepsilon>0, \delta>0$.
\end{example}}

It is well-known that  the limiting subdifferential can be used to characterize Lipschitz properties of functions. 
\begin{lemma}\cite[Theorem 9.13]{RW}) \label{Lemma2.1} Let $f:\mathbb R^n\rightarrow \overline{\mathbb R}$ with $f(\bar x) $  finite. Suppose that $f$ is locally l.s.c. at $\bar x$.
Then the following conditions are equivalent:
\begin{itemize}
\item[{\rm (a)}]
$f$ is Lipschitz continuous at $\bar x$,
\item[{\rm (b)}]$\partial^\infty f(\bar x)=\{0\}$,
\item[{\rm (c)}] the set-valued map  $\partial f$ is locally bounded at $\bar x$, i.e.,  there exists $\varepsilon>0$ such that 
$\bigcup_{x\in  \mathbb{B}_\varepsilon(\bar x)} \partial f( x) $ is bounded.
\end{itemize}
Moreover, when these conditions hold, $\partial f(\bar x)$ is nonempty and compact and one has
$$ {\rm lip} f(\bar x)=\max_{v\in \partial f(\bar x)}\|v\|:= \max |\partial f(\bar x)|,$$
where ${\rm lip} f(\bar x)$ is the Lipschitz modulus of $f$ at $\bar x$ defined by
\begin{equation}
{\rm lip}f(\bar x):=\displaystyle \limsup_{x',x\rightarrow \bar x,\atop x'\not =x}\frac{\|f(x')-f(x)\|}{\|x-x'\|}.
\label{Lipsnew}
\end{equation} 
\end{lemma}
We now prove that the directional limiting subdifferential can also be used to characterize the directional Lipschitz continuity. Note that when the direction $u\not =0$, we need the continuity of $f$ in direction $u$.
\begin{prop}\label{liplemmanew}
	Let $f:\mathbb R^n\rightarrow \overline{\mathbb R}$ with $f(\bar x) $  finite. Assume that $f$ is locally l.s.c. at $\bar x$ in direction $u$ and $f$ is continuous at $\bar x$ in direction $u$.
{Then $f(x)$ is Lipschitz continuous around $\bar x$ in direction $u$, if and only if  $\partial f$ is uniformly bounded on $(\bar x+{\mathcal V}_{ \varepsilon,\delta}(u))\backslash\{\bar x\}$, i.e.,
\begin{equation} \max_{x\in (\bar x+{\mathcal V}_{\varepsilon,\delta}(u))\backslash\{\bar x\}} |\partial f(x)|:=\max_{v\in \partial f( x) \atop{x\in (\bar x+{\mathcal V}_{\varepsilon,\delta}(u))\backslash\{\bar x\}}} \|v\| \leq M, \label{uniformbnew} \end{equation} 
for some positive scalars $\varepsilon,\delta, M$.} 
\end{prop}
\beginproof 
Since $f$ is locally l.s.c. at $\bar x$ in direction $u$, there exist positive scalars $\varepsilon,\delta_1$ such that epi$f\cap \overline{\bar x+{\mathcal V}_{\varepsilon,\delta_1}(u)}\times \overline{\mathbb B}_\varepsilon(f(\bar x)))$ is closed.  Take $\delta<\delta_1$ and consider points  $x$ lying in  $(\bar x+{\mathcal V}_{\varepsilon,\delta}(u))\backslash\{\bar x\}$.  Then $x$ is an  interior point of  $(\bar x+{\mathcal V}_{\varepsilon,\delta_1}(u))\backslash\{\bar x\}$, and we can choose a small enough $\bar{\varepsilon}>0$ such that epi$f\cap \overline{\mathbb{B}}_{\bar{\varepsilon}} (x, f(x))$ is closed which means that $f$ is locally l.s.c. at $x$. Moreover since $f$ is continuous at $\bar x$ in direction $u$, without loss of generality we may assume that $f$ is finite on $(\bar x+{\mathcal V}_{\varepsilon,\delta}(u))\backslash\{\bar x\}$. Hence without loss of generality, we may assume that $f$ is finite and locally l.s.c. at each point lying in $(\bar x+{\mathcal V}_{\varepsilon,\delta}(u))\backslash\{\bar x\}$ for some positive scalars $\varepsilon,\delta$.

Suppose that $f$ is Lipschitz continuous at $\bar x$ in direction $u$, then
at each $ \tilde x\in (\bar x+{\mathcal V}_{\varepsilon,\delta}(u))\backslash\{\bar x\}$, $f$ is Lipschitz continuous and $f(\tilde x)$ is finite.  By   Lemma \ref{Lemma2.1},  $\partial f$ is locally bounded at each $\tilde x$ with ${\rm lip} f(\tilde x)\leq M$ for some $M>0$. By (\ref{Lipsnew}), we obtain the uniform boundedness of  $\partial f$ on $(\bar x+{\mathcal V}_{ \varepsilon,\delta}(u))\backslash\{\bar x\}$ for sufficiently small positive scalars $\varepsilon,\delta$. 
 Conversely suppose that $\partial f$ is uniformly bounded on $(\bar x+{\mathcal V}_{ \varepsilon,\delta}(u))\backslash\{\bar x\}$ and (\ref{uniformbnew}) holds. Then $\partial f$ is locally bounded at each  point  $\tilde x\in (\bar x+{\mathcal V}_{\varepsilon,\delta}(u))\backslash\{\bar x\}$. By  Lemma \ref{Lemma2.1}, local boundedness of $\partial f$ at $\tilde x$ {is equivalent to the Lipschitz continuity of $f$ at $\tilde x$ and in this case ${\rm lip} f(\tilde x)\leq M$. It follows that there exists $\tilde\varepsilon>0$ such that
 {
  \begin{equation*}
 	|f(x')-f(x'')|\leq M \|x'-x''\|, \quad \forall x',x''\in \mathbb B_{\tilde\varepsilon}(\tilde x).
 \end{equation*}	
Then for any $a,b\in(\bar x+{\mathcal V}_{\varepsilon,\delta}(u))\backslash\{\bar x\}$, by the convexity of set $\bar x+{\mathcal V}_{\varepsilon,\delta}(u)$ the segment $[a,b]$ joining  $a$ and $b$ is contained in $\bar x+{\mathcal V}_{\varepsilon,\delta}(u)$. Since $[a,b]$ is compact, by the Heine-Borel theorem, one can find finitely many points, say $x^2,\ldots,x^n$ from the segment $[a,b]$, such that 
\[
\|a-b\|=\|a-x^2\|+\|x^2-x^3\|+\ldots+\|x^n-b\|
\]
and $|f(x^i)-f(x^{i+1})|\leq M \|x^i-x^{i+1}\|$ for $i=1,\ldots,n$ where $x^1:=a, x^{n+1}:=b$. It follows that 
\begin{align*}
|f(a)-f(b)|&\leq|f(x^1)-f(x^2)|+\ldots+|f(x^n)-f(x^{n+1})|\\
&\leq M\|x^1-x^2\|+\ldots+M\|x^n-x^{n+1}\|\\
&=M\|a-b\|,
\end{align*}
from which we have
\begin{align}\label{lip}
	|f(a)-f(b)|\leq M \|a-b\|, \quad \forall a,b\in (\bar x+{\mathcal V}_{\varepsilon,\delta}(u))\backslash\{\bar x\}.
\end{align}
In (\ref{lip}) let $a$ approach $\bar x$ in direction $u$. Then by the directional continuity of $f(x)$ at $\bar x$ in direction $u$, we obtain
\[
|f(\bar x)-f(b)|\leq M\|\bar x-b\|, \quad \forall b\in (\bar x+{\mathcal V}_{\varepsilon,\delta}(u))\backslash\{\bar x\}.
\]The proof is complete.
}
\endproof

Using Proposition \ref{liplemmanew}, we can show that the directional singular subdifferential can be used to characterize the directional Lipschitz continuity of a directionally l.s.c. function.
\begin{prop}\label{suflip}
			Let $f:\mathbb R^n\rightarrow \overline{\mathbb R}$ with $f(\bar x)$ finite be {locally} l.s.c. and 
			 continuous at $\bar x$ in direction $u\not =0$.   Then $f(x)$ is Lipschitz continuous around $\bar x$ in direction $u$ if and only if $\partial^\infty f(\bar x;u)=\{0\}$.
\end{prop}
\beginproof
If $f(x)$ is Lipschitz around $\bar x$ in direction $u$, then by \cite[Corollary 5.9]{Long}, we have $\partial^\infty f(\bar x;u)=\{0\}$. 
Suppose that $\partial^\infty f(\bar x;u)=\{0\}$. Then by definition,  there  exist positive scalars $\varepsilon,\delta$ such that $\partial f$ is uniformly bounded on $(\bar x+{\mathcal V}_{\varepsilon,\delta}(u))\backslash\{\bar x\}$,  by Proposition \ref{liplemmanew}, $f(x)$ is Lipschitz continuous around $\bar x$ in direction $u$. 
\endproof

{The following proposition recalls some useful calculus rules of directional and nondirectional limiting subdifferentials. The readers are also refered to \cite{RW} and \cite{BHA} for comprehensive studies of calculus properties of these two kinds of subdifferentials.}
\begin{prop}[Calculus Rules]{(see e.g., \cite{BY, BP2022, RW, Aub2,Long})} \label{Calculus}
\begin{itemize}
\item[{\rm (1)}] Let $f:\mathbb{R}^n \rightarrow\mathbb{R}$ be Lipschitz around $\bar{x}$ in direction $u$ and $g:\mathbb{R}^n\rightarrow 
\overline{\mathbb{R}}$ be l.s.c. in direction $u$ with $\bar x\in {\rm dom} g$. Let $\alpha,\beta$ be nonnegative scalars.  Then $$\partial (\alpha f+ \beta g)(\bar{x};u)\subseteq \alpha\partial f(\bar x;u)+\beta\partial g(\bar{x};u).$$
\item[{\rm (2)}]  Let $\phi:\mathbb{R}^n\rightarrow \mathbb{R}^m$ be Lipschitz near $\bar{x}$ and $f:\mathbb{R}^m
\rightarrow\mathbb{R}$ be Lipschitz near $\phi(\bar{x})$. Then
    \begin{equation*}
    \partial(f\circ\phi)(\bar{x})\subseteq\bigcup_{\xi\in\partial f(\phi{(\bar{x})})}\partial\langle\xi,\phi\rangle(\bar{x}).
    \end{equation*}
\item[{\rm (3)}] Let $f:\mathbb{R}^n\rightarrow\mathbb{R}$ be Lipschitz around $\bar{x}$ and $\Omega$ be a closed subset of $\mathbb{R}^n$.
If $\bar{x}$ is a local minimizer of $f$ over $\Omega$, then $0\in\partial f(\bar{x})+N_\Omega(\bar{x})$.
\end{itemize}
\end{prop}

Based on Definition \ref{dn} we can give the definition of directional metric subregularity/regularity.
\begin{defn}[Directional Metric Subregularity/Regularity] \cite[Definition 1]{Gfr13}
	Let $G:\mathbb{R}^n\rightrightarrows\mathbb{R}^m$ be a multifunction given by $G(x):=\Gamma-\phi(x)$, where $\phi:\mathbb R^n\rightarrow\mathbb R^m$ and $\Gamma\subseteq \mathbb R^m$ is closed. Further suppose that  $\bar{y} \in \Gamma-\phi(\bar{x})$, $u\in\mathbb{R}^n$ and $v\in \mathbb{R}^m$.
	\begin{itemize}
\item[1.] We say that  the set-valued map
$G$ is metrically subregular  at $(\bar{x},\bar{y})$ in direction $u$ or the system $\phi(x)\in \Gamma$ is calm at $(\bar{x},\bar{y})$ in direction $u$, 
if  there are positive reals $\varepsilon>0,\delta>0,$ and $\kappa'>0$ such that
		\begin{equation*}
		d(x,\phi^{-1}(\Gamma-\bar{y}))\leq\kappa' d(\phi(x)+\bar y,\Gamma),
		\end{equation*}
		where $\phi^{-1}(\Gamma-\bar y):=\{x| \bar y\in \Gamma- \phi(x)\}$
		holds for all $x\in\bar{x}+{\mathcal V}_{\varepsilon,\delta}(u)$. When $u=0$ in the above definition, we say that the set-valued map $G$ is  metrically subregular   at $(\bar{x},\bar{y})$ or the system $\phi(x)\in \Gamma$ is calm at $\bar x$.
\item[2.] 	
We say that the set-valued map
$G$ is metrically regular  at $(\bar{x},\bar{y})$ in direction $(u,v)$,
if there are positive reals $\varepsilon>0,\delta>0,$ and $\kappa'>0$ such that
		\begin{equation*}
		d(x,\phi^{-1}(\Gamma-y))\leq\kappa' d(\phi(x)+ y,\Gamma),
		\end{equation*}
		holds for all $(x,y) \in (\bar{x},\bar{y})+
		{\mathcal V}_{\varepsilon,\delta}(u,v)$ with $d(x,\phi^{-1}(\Gamma-y))
		\leq \delta \|(x,y)-(\bar{x},\bar{y})\|$ if $\|(u,v)\|\not =0$. When $(u,v)=(0,0)$ in the above definition, we say that the set-valued map $G$ is  metrically regular at $(\bar{x},\bar{y})$.
		\end{itemize}	
\end{defn}
Consider the case where $\phi$ is  Lipschitz around $\bar x$. It is known that $G(x)=\Gamma-\phi(x)$ is metrically regular at $(\bar x, 0)$  if and only if the Mordukhovich's criterion/no nonzero abnormal multiplier constraint qualification (NNAMCQ)  holds (see e.g., \cite{RW}), i.e.
$$ 0\in\partial \langle \phi,\lambda\rangle (\bar x) \quad \lambda \in N_\Gamma(\phi(\bar x)) \Longrightarrow \lambda =0,$$
 which is equivalent to MFCQ for the case of smooth nonlinear programs. 
Recently the so-called first-order sufficient condition for directional metric subregularity (FOSCMS) was introduced in \cite{GK} when $\phi$ is smooth and extended to nonsmooth case in  \cite{BHA}, which says that if for all $w\in D\phi(\bar x)(u)$ with $w\in T_\Gamma(\phi(\bar x))$ one has implication
$$ 0\in\partial \langle\phi,\lambda\rangle(\bar x;u)=0,\quad \lambda \in N_\Gamma(\phi(\bar x); w) \Longrightarrow \lambda=0,$$
then $G(x)=\Gamma-\phi(x)$ is metrically subregular at $(\bar x,0)$ in direction $u$.  
  Another useful criterion for the metric subregularity which requires that  $\phi$ is affine and $\Gamma$ is the union of finitely many convex polyhedral sets,  is based on Robinson's multifunction theory \cite{Robinson}. There are also other weaker criteria for metric subregularity such as the quasi-/pseudo-normality (\cite{GYZ}) and the directional  quasi-/pseudo-normality \cite{BYZ,BP2022A}. These criteria are weaker but are sequential, not point-based.

The following first-order sufficient condition for metric regularity that we will need to use in this paper is a slightly stronger condition than \cite[Theorem 6.1(2.)]{BHA}. For completeness we provide the proof. We {denote by $\tilde\partial f(\bar x;(u,\xi))$, the directional limiting subdifferential of function $f$ at $\bar x$ in direction $(u,\xi)\in T_{\mbox{gph}f}(\bar x,f(\bar x))$}  (see  \cite{BHA})}, and  use this notation only in Propositions \ref{dMR}, \ref{chainrule} and Remark \ref{remark3.2}.
{\begin{prop}[First-order Sufficient Condition for  Metric Regularity] \label{dMR}
Let $\phi:\mathbb R^n\rightarrow\mathbb R^p$ with $\Gamma\subseteq \mathbb R^p$ closed. Suppose that  $\phi(\bar x) \in \Gamma$ and  $\phi$ is locally Lipschitz continuous around $\bar x$. Then 
the  set-valued map $G(x):={\Gamma-\phi(x)}$ is metrically regular  at $(\bar x,0)$ in direction $(u,v)$ provided for all 
$w\in D\phi(\bar x)(u)$ with $w+v\in  T_\Gamma(\phi(\bar x))$ one has the implication
\begin{equation}\label{aFOSCMS}
	0\in\partial \langle\lambda,\phi\rangle(\bar x;u)=0,\quad \lambda \in N_\Gamma(\phi(\bar x); w+v) \Longrightarrow \lambda=0.
\end{equation}
\end{prop}
\beginproof
By Benko et al. \cite[Theorem 6.1(2.)]{BHA}, the metric regularity for the set-valued map $G(x)=\Gamma-\phi(x)$ holds at $(\bar x,0)$ in direction $(u,v)$ provided  \begin{eqnarray}\label{FOSCMS}
	&& 0\in D^*\phi(\bar x;(u,w))(\lambda),\ \lambda \in N_\Gamma(\phi(\bar x); w+v) \nonumber \\
	&& \implies \lambda=0,\ \mbox{for all}\ w\in D\phi(\bar x)(u)\ \mbox{with}\ w+v\in T_\Gamma(\phi(\bar x)),
\end{eqnarray}where $D^*\phi(\bar x;(u,w))$ is the limiting coderivative of $\phi$ at $\bar x$ in direction $(u,w)\in T_{\mbox{gph}\phi}(\bar x,\phi(\bar x))$ defined in  \cite{BHA}.
By \cite[Proposition 5.1]{BHA}, since $\phi(x)$ is  Lipschitz continuous near $\bar x$, one has
\[
D^*\phi(\bar x;(u,w))(\lambda)=\tilde\partial\langle\lambda,\phi\rangle(\bar x;(u,\lambda^Tw)).
\]
Then condition (\ref{FOSCMS}) is equivalent to 
\begin{eqnarray*}
&& 0\in \tilde\partial\langle\lambda,\phi\rangle(\bar x;(u,\lambda^Tw)),\ \lambda \in N_\Gamma(\phi(\bar x);  w+v) \\
&&  \Longrightarrow \lambda=0,\ \mbox{for all}\ w\in D\phi(\bar x)(u)\ \mbox{with}\ w+v\in T_\Gamma(\phi(\bar x)),
\end{eqnarray*}
which is weaker than condition
\begin{align*}
	&0\in \bigcup_{{w'}\in D\phi(\bar x)(u)}\tilde\partial\langle\lambda,\phi\rangle(\bar x;(u,\lambda^T {w'})),\ \lambda \in N_\Gamma(\phi(\bar x);  w+v) \\
	 &\Longrightarrow \lambda=0,\ \mbox{for all}\ w\in D\phi(\bar x)(u)\ \mbox{with}\ w+v\in T_\Gamma(\phi(\bar x)).
\end{align*} 
By \cite[Corollary 4.1]{BHA}, since $\phi(x)$ is  Lipschitz continuous near $\bar x$, one has
\[
\partial\langle\lambda,\phi\rangle(\bar x;u)=\bigcup_{v\in D\langle\lambda,\phi\rangle(\bar x)(u)}\tilde\partial\langle\lambda,\phi\rangle(\bar x;(u,v))\bigcup_{w'\in D\phi(\bar x)(u)}\tilde\partial\langle\lambda,\phi\rangle(\bar x;(u,\lambda^Tw')),
\]{where the second equality can be easily obtained from the Lipschitzness of $\phi(x)$ and Definition \ref{derivatives}.}
Consequently, condition (\ref{aFOSCMS}) implies condition (\ref{FOSCMS}). The proof is complete.
\endproof}

Like its nondirectional version, the directional metric subregularity can be used to derive calculus rules of directional limiting subdifferentials. The following chain rule will be used in this paper.

{\begin{prop}\cite[Theorem 4.1, Corollary 4.1, Proposition 5.1]{BHA}\label{chainrule}
	Let $\phi:\mathbb R^n\rightarrow\mathbb R^p$ be directionally Lipschitz continuous at $\bar x$ in direction $u$, and $\Gamma\subseteq \mathbb R^p$ be closed. Suppose $\phi(\bar x) \in \Gamma$, $D\phi(\bar x)(u)\cap T_\Gamma(\phi(\bar x))\neq\emptyset$ and the metric subregularity for the set-valued map  $G(x):=\Gamma-\phi(x)$ holds at $(\bar x,0)\in {\rm gph}G$ in direction $u$. Then 
	\begin{equation}\label{chain1}
	\partial (\delta_\Gamma\circ \phi)(\bar x;u)\subseteq\bigcup_{d\in D\phi(\bar x)(u)\cap T_\Gamma(\phi(\bar x))}\{\partial\langle\zeta,\phi\rangle(\bar x;u)|\zeta\in N_\Gamma(\phi(\bar x);d)\}.
	\end{equation}
\end{prop}
\beginproof
Since
\begin{align}\label{indic}
	D\delta_\Gamma(\phi(\bar x))(d)=\left\{
	\begin{array}{ll}
		\{0\},\ \mbox{if}\ d\in T_\Gamma(\phi(\bar x)),\\
		\emptyset,\ \mbox{otherwise},
	\end{array}
	\right.
\end{align} and by \cite[Corollary 4.1]{BHA},
\begin{equation}\label{equiv}
\partial (\delta_\Gamma\circ \phi)(\bar x;u)=\tilde\partial (\delta_\Gamma\circ \phi)(\bar x;(u,0)).
\end{equation}
 Since $\phi(x)$ is Lipschitz continuous near $\bar x$, by \cite[Theorem 4.1]{BHA} and (\ref{indic}), one has the directional limiting subdifferential chain rule
\begin{equation}\label{chain}
	\tilde\partial (\delta_\Gamma\circ \phi)(\bar x;(u,0))\subseteq \bigcup_{d \in D\phi(\bar x)(u)\cap T_\Gamma(\phi(\bar x))}D^*\phi(\bar x;(u,d))\tilde\partial\delta_\Gamma(\phi(\bar x);(d,0)).
\end{equation}
By (\ref{indic}) and since $\tilde\partial\delta_\Gamma(\phi(\bar x);(d,0))=\partial\delta_\Gamma(\phi(\bar x);d)=N_\Gamma(\phi(\bar x);d)$ by \cite[Formula (2)]{BHA}, condition (\ref{chain}) is equivalent to 
\begin{equation*}
	\tilde\partial (\delta_\Gamma\circ \phi)(\bar x;(u,0))\subseteq \bigcup_{d\in D\phi(\bar x)(u)\cap T_\Gamma(\phi(\bar x))}D^*\phi(\bar x;(u,d))N_\Gamma(\phi(\bar x);d).
\end{equation*}
Then similar to the proof of Proposition \ref{dMR}, by \cite[Corollary 4.1]{BHA} and (\ref{equiv}), one can obtain (\ref{chain1}). The proof is complete.
\endproof}

\section{Directional subdifferentials of the value function}
In this section we give upper estimates for the directional limiting and the singular subdifferentials of the value function. From these estimates we obtain  sufficient conditions for directional Lipschitz continuity of the value function.

{The first desired property one wishes to have for the value function is the lower semicontinuity. 
A very weak sufficient condition for the lower semicontinuity was introduced by Clarke in \cite[Hypothesis 6.5.1]{Clarke}. This condition was referred as the restricted inf-compactness condition  in \cite{GLYZ}. The following   directional version of the restricted inf-compactness condition was introduced in \cite[Definition 4.1]{BY}.}
\begin{defn}[Directional Restricted Inf-compactness]\cite[Definition 4.1]{BY}\label{reinf}
	We say that the restricted inf-compactness holds at $\bar x$ in direction $u$ if $V(\bar x)$ is finite and there exist a compact set  $\Omega_u\subseteq \mathbb{R}^n$, and positive numbers $\varepsilon>0, \delta>0$ such that for all $ x\in \bar x+{\mathcal V}_{\varepsilon,\delta}(u)$ with $V(x)<V(\bar x)+\varepsilon$, one always has $S(x)\cap\Omega_u\neq\emptyset$.
\end{defn}
Obviously, if the restricted inf-compactness holds at $\bar x$ in {direction u=0}, then the {classical} restricted inf-compactness  holds.
{The restricted inf-compactness at $\bar x$ in direction $u$ is weaker than the inf-compactness condition at $\bar x$ in direction $u$,  by which we mean that there exist 
$\alpha, \varepsilon>0,\delta >0$ and a bounded set $\Omega_u$ such that  $\alpha >V(\bar x)$ and  \[
\{y\in\mathbb R^m|P( x,y)\in\Gamma, f( x,y)\leq \alpha, x \in \bar x+  {\mathcal V}_{\varepsilon,\delta}(u)\}	\subseteq \Omega_u.
\] 
Clearly, if the feasible region is uniformly bounded around $\bar x$ in direction $u$, by which we mean the set
{$$ \bigcup_{x\in \bar x+  {\mathcal V}_{\varepsilon,\delta}(u)} {\mathcal F}(x)=\{y\in\mathbb R^m|P(x,y)\in\Gamma, x \in \bar x+  {\mathcal V}_{\varepsilon,\delta}(u)\}$$}
is bounded, then the inf-compactness condition  at $\bar x$ in direction $u$ holds.}

It is known that under the restricted inf-compactness condition, the value function is locally lower semicontinuous (see Clarke \cite[Page 246]{Clarke}). The following proposition gives a directional version of this result.
\begin{prop} \label{directional-LSC} Suppose that the restricted inf-compactness holds at $\bar x$ in direction $u$. Then $V(x)$ is locally l.s.c. at $\bar x$ along direction $u$. 
\end{prop}
\beginproof 
By Definition \ref{locallsc}, we need to prove that there exist some positive scalars $\varepsilon,\delta$ such that  the set
\begin{equation*}
{\rm epi}V\cap \overline{\bar x+{\mathcal V}_{\varepsilon,\delta}(u)}\times \overline{V(\bar x)+\mathbb{B}}_\varepsilon
\end{equation*}
is closed.

Let $\varepsilon, \delta>0$ be given as in Definition \ref{reinf}. 
Let $(x^k, r^k )\in {\rm epi}V\cap \overline{\bar x+{\mathcal V}_{\varepsilon,\delta}(u)}\times \overline{V(\bar x)+\mathbb{B}}_\varepsilon$ and suppose that $x^k\rightarrow \tilde{x}$ and $r^k\rightarrow \tilde{r}$. By the closedness of the set $\overline{\bar x+{\mathcal V}_{\varepsilon,\delta}(u)}\times \overline{V(\bar x)+\mathbb{B}}_\varepsilon$, $(\tilde{x},\tilde{r}) \in \overline{\bar x+{\mathcal V}_{\varepsilon,\delta}(u)}\times \overline{V(\bar x)+\mathbb{B}}_\varepsilon.$ So we only need to show that $V(\tilde x)\leq \tilde{r}$.
Since $V(x^k)\leq r^k$ by definition, we have $ V(x^k) \leq r^k \leq V(\bar x)+ \varepsilon$. Then by Definition \ref{reinf}, for sufficiently large $k$, there exists $y^k\in S(x^k)\cap\Omega_u$. By the compactness of $\Omega_u$, without loss of generality, we may assume that   $y^k\rightarrow\tilde y$ for some $\tilde y\in \Omega_u$. Since $P(x^k,y^k)\in \Gamma$, by the continuity of $P(x,y)$ and the closedness of set $\Gamma$, we have $P(\tilde x,\tilde y)\in \Gamma$. 
It follows that 
$$V(\tilde x) \leq f(\tilde x, \tilde y) =\lim_{k\rightarrow \infty} f(x^k,y^k) =\lim_{k\rightarrow \infty} V(x^k) \leq \lim_{k\rightarrow \infty} r^k =\tilde r,$$
and hence the proof is complete.
\endproof

Combing the results in Propositions \ref{suflip} and \ref{directional-LSC}, we obtain the following sufficient condition for directional Lipschitz continuity of the value function by using the directional singular subdifferential of the value function.
\begin{prop}\label{VLip}
Suppose the restricted inf-compactness condition holds at $\bar x$ in direction $u$ and $V(x)$ is continuous at $\bar x$ in direction $u$ if $u\not =0$.  Then $V(x)$ is Lipschitz continuous around $\bar x$ in direction $u$ if and only if $\partial^\infty V(\bar x;u)=\{0\}$.
\end{prop}

In the following definition, we give a subset of the solution set based on a direction. This set can be used to provide a tighter upper estimates for the directional subdifferentials than the solution set as in Theorems \ref{estimates} and \ref{isestimates}.
\begin{defn}[Directional Solution]\cite[Definition 4.5]{BY}
	The optimal solution in direction $u$ is defined by
	$$S(\bar x;u)=\{y\in S(\bar x)| \exists t_k\downarrow 0, u^k\rightarrow u, y^k\rightarrow y, y^k\in S(\bar x+t_k u^k)\}.$$  If $\bar y\in S(\bar x;u)$, we say that $\bar y$ is upper stable in direction $u$ in the sense of Janin (see \cite[Definition 3.4]{Janin}).
\end{defn}

The following directional version of inner semicontinuity can refine our results.
\begin{defn}[Directional Inner Semicontinuity]\cite{BHA} Given $(\bar x,\bar y)\in \mbox{gph} S$, we say that  the optimal solution map $S(x)$ is inner semicontinuous at $(\bar x,\bar y)$ in direction $u$, if for any sequences $x^k\xrightarrow{u}\bar x$, there exists a sequence $y^k\in S(\bar x+t_ku^k)$ converging to $\bar 
	y$. 
\end{defn}
By definition, if $\exists \bar y\in S(\bar x)$ such that $S(x)$ is inner semicontinuous  at $(\bar x,\bar y)$ in direction $u$, then the restricted inf-compactness holds at $\bar x$ in direction $u$. Note that if $S(x)$ is inner semicontinuous at $(\bar x,\bar y)$ in {direction u=0}, then $S(x)$ is inner semicontinuous at $(\bar x,\bar y)$ in the sense of \cite[Definition 1.63]{Aub2}.



In the following concepts, the inner semicontinuity properties are strengthened by controlling the rate of  convergence
$y^k\rightarrow \bar y$. \begin{defn}[Directional Inner Calmness/Calmness*](Benko et al. \cite{B2021}) The optimal solution map $S(x)$ is said to be
	\begin{itemize}
		\item[{\rm (i)}] inner calm* at $\bar x$ in direction $u$ if there exists $\kappa>0$ such that for every sequence $x^k\xrightarrow{u}\bar x$, 
		there exist a subsequence $K$ of the set of nonnegative integers $\mathbb{N}$, together with a sequence  $ y^k\in S(x^k)$ for $k\in K$ and $\bar y\in\mathbb R^m$ such that 
		\begin{equation}\label{ic}
		\| y^k-\bar y\|\leq\kappa\| x^k-\bar x\|.
		\end{equation} 
		\item[{\rm (ii)}] inner calm at $(\bar x,\bar y)\in \mbox{gph} S$ in direction $u$ if there exist $\kappa>0$ such that for every sequence $x^k\xrightarrow{u}\bar x$ there exists a sequence $y^k$ satisfying  $y^k\in S(x^k)$ and (\ref{ic}) 
		 for sufficiently large $k$, {or equivalently there exist $\kappa>0,\varepsilon>0,\delta>0$ such that
		$$\bar y \in S(x) +\kappa \|x-\bar x\|\mathbb{B} \qquad \forall x\in \bar x+{\mathcal V}_{\varepsilon,\delta}(u).$$}
	\end{itemize}
\end{defn}

\begin{figure}[h]
	\centering
	\includegraphics[width=15cm]{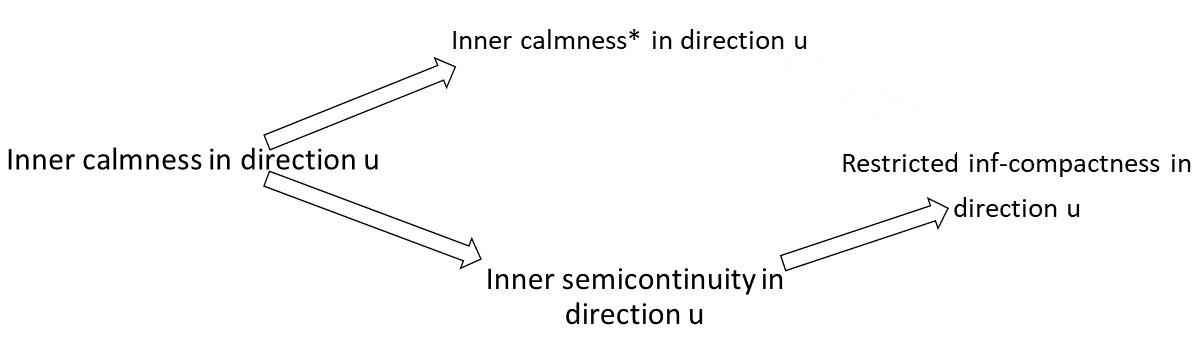}
	\caption{Implications}
	\label{rel}
\end{figure}

{ Similar  to \cite{B2021}, one can easily obtain the following implications in Figure \ref{rel}.} Note that we can not derive the relationship between the inner calmness* and the restricted inf compactness  of optimal solution map $S(x)$ at $\bar x$ in direction $u$ since in the definition of  the inner calmness* condition, $\bar y$ does not necessarily belong to $S(\bar x)$.

Sufficient conditions for the inner semicontinuity of optimal solution map $S(x)$ of parametric optimization problems can be found in \cite[Theorem 5.9]{RW} and \cite[Remark 3.2]{DDM2007}. 
Now we list some  sufficient conditions for the (directional) inner calmness*/semicontinuity/calmness.  
\begin{itemize}
	\item[(a)]  {If $S(x)$ is a continuous single-valued map around $\bar x$, then $S(x)$ is inner semicontinuous at $(\bar x, S(\bar x))$. For example, for parametric nonlinear programs, if the linear independence constraint qualification, the sufficient second-order optimality condition and the strict complementarity condition hold at $(\bar x,\bar y,\bar\mu)$ with $\bar y\in S(\bar x)$ and $\bar\mu$ being the associated KKT multiplier, then $S(x)$ is single-valued and continuous near $\bar x$, see \cite[Lemma 2.1]{DZ}. } 
	\item[(b)] If $S(x)$ is Lipschitz-like (or has the Aubin property) around $(\bar x,\bar y)\in \mbox{gph} S$ in direction $u$, i.e., there exist $\kappa>0,\varepsilon>0,\delta>0$  and $V$, a neighborhood of $\bar y$ such that 
	\begin{equation*}
			S(x)\cap V\subseteq  S(x') +\kappa \|x-x'\|\mathbb{B} \qquad \forall x, x'\in \bar x+{\mathcal V}_{\varepsilon,\delta}(u),
		\end{equation*} then $S(x)$ is inner calm at  $(\bar x,\bar y)$ in direction $u$. In \cite{Aub2,LM2004}, various sufficient conditions for the Aubin property of solution maps are established. 
	\item[(c)] If the graph of $S(x)$ is convex and $\bar x\in\mbox{int}(\mbox{dom} S)$, then $S(x)$ is inner semicontinuous at $\bar x$ \cite[Theorem 5.9(b)]{RW}.
	\item[(d)]  Recall that a set-valued map is called polyhedral if its graph is a union of finitely many convex polyhedral sets. It was shown in \cite[Theorem 3.4]{B2021} that polyhedral maps enjoy the inner calmness* property. It follows that the solution set of a linear program with additive perturbation $S(x):=\displaystyle \arg\min_y \{c^T y| Ay\leq x\}$
	where $c\in \mathbb{R}^m, A\in \mathbb{R}^{n\times m}$
	 is inner calm*  at any $\bar x$. Also  $S(x)$ is inner semicontinuous at any point in the graph by 
	 \cite[Theorem 4.3.5]{Bank}.
\end{itemize}


In Proposition \ref{directional-LSC} we have shown that the directional restricted inf-compactness condition implies the directional local lower semicontinuity of $V(x)$. {Although the directional inner calmness* and the directional restricted inf-compactness condition are not comparable, following a similar proof as Proposition \ref{directional-LSC}, one can show that the directional inner calmness* implies the directional local lower semicontinuity. Hence, by Proposition \ref{suflip}, we also have the following sufficient condition for the Lipschitz continuity of the value function.
\begin{prop}\label{Vlip2}
	Suppose that the inner calmness* holds at $\bar x$ in direction $u$ and $V(x)$ is continuous at $\bar x$ in direction $u$ if $u\not =0$.  Then $V(x)$ is Lipschitz continuous around $\bar x$ in direction $u$ if and only if $\partial^\infty V(\bar x;u)=\{0\}$.
\end{prop} In the next proposition, we show that  under the directional inner semicontinuity/calmness which are stronger than the directional restricted inf-compactness condition, $V(x)$ is directionally continuous.
\begin{prop}
\label{directional-conti}
$V(x)$ is continuous at $\bar x$ in direction $u$ under one of the following assumptions.
	\begin{itemize} 
		\item[{\rm (1)}]  $S(x)$ is inner semicontinuous at $(\bar x,\bar y)\in \mbox{gph} S$ in direction $u$. 
		\item[{\rm (2)}]  $S(x)$ is inner calm at $(\bar x,\bar y)\in \mbox{gph} S$ in direction $u$. 
	\end{itemize}
\end{prop}
\beginproof
{(1) Consider any sequence $x^k\xrightarrow{u}\bar x$. Since $S(x)$ is inner semicontinuous at $(\bar x,\bar y)\in \mbox{gph} S$ in direction $u$, there exists a sequence $y^k\in S(x^k)$ such that $y^k\rightarrow\bar y$.
This implies that $\lim_k V(x^k)=\lim_kf(x^k,y^k)=f(\bar x,\bar y)=V(\bar x)$.
Finally, by the choice of the sequence $\{x^k\}$, one has $\lim_{x\xrightarrow{u}\bar x}V(x)=V(\bar x)$. Hence, $V(x)$ is continuous at $\bar x$ in direction $u$.
Since the directional inner calmness implies the  directional inner semicontinuity, (2) holds immediately from (1).
\endproof

Denote the  the linearization cone of {\rm gph}${\mathcal F}$ at $(x,y)\in {\rm gph}{\mathcal F} $ as
$$\mathbb{L}(x,y):=\{(u,v)|DP(x,y)(u,v)\cap T_\Gamma(P(x,y))\neq\emptyset\},$$
and its $y$-projection in direction $u$ as
$$\mathbb{L}(x,y;u):=\{v\mid DP(x,y)(u,v)\cap T_\Gamma(P(x,y))\neq\emptyset\}.$$
Note that when $u=0$, $\mathbb{L}(x,y;u)$ becomes the  linearization cone at $y$ for problem $(P_x)$
$$\mathbb{L}(x,y;0)=\{v|DP(x,y)(0,v)\cap T_\Gamma(P(x,y))\neq\emptyset\}.$$

The following proposition gives the relationship between the tangent cone of the feasible region and its linearization cone.
\begin{prop}\label{Abadie}
	Let $\bar x$ be a feasible point of the system $\phi(x) \in \Gamma$, where $\phi$ is Lipschitz around $\bar x$ and $\Gamma$ is closed.  One  has \begin{equation}\label{linear}
		T_{\phi^{-1}(\Gamma)}(\bar x)\subseteq \mathbb{L}(\bar x):=\{u|D\phi(\bar x)(u)\cap T_\Gamma(\phi(\bar x))\neq\emptyset\}.
	\end{equation} Furthermore, if the metric subregularity of the set-valued map $\Psi(x):=\Gamma-\phi(x)$ holds at $(\bar x, 0)$, and if either $\phi(x)$ is directionally differentiable at $\bar x$ in direction $u$ or the set $\Gamma$ is geometrically derivable at {$\phi(\bar x)$}, then the equality in (\ref{linear}) holds.
\end{prop}
\beginproof
Take any direction $u\in T_{\phi^{-1}(\Gamma)}(\bar x)$. Then there exist sequences $t_k\downarrow0, u^k\rightarrow u$ such that $\bar x
+t_ku^k\in {\phi^{-1}(\Gamma)}$. Equivalently, $\phi(\bar x+t_ku^k)\in \Gamma$. Then by the Lipschitz continuity of $\phi(x)$, passing to a subsequence if necessary, there exists a vector 
\[
d:=\lim_{k\rightarrow \infty}\frac{\phi(\bar x+t_ku^k)-\phi(\bar x)}{t_k}.
\]Obviously, $d\in D\phi(\bar x)(u)\cap T_\Gamma(\phi(\bar x))$. Hence, $
u\in\mathbb L(\bar x)$.
The opposite inclusion follows by combining \cite[Proposition 4.1]{GO2016} and \cite[Proposition 2.2]{BYZ}.
\endproof

{The upper/lower Dini derivative of the value function $V(x)$ is employed in the following results. The reader is referred to \cite{BY} and \cite{GLYZ, LY, LYErru} for formulae/estimates of Dini/directional derivative of the value functions of parametric nonlinear programs/mathematical programs with equilibrium constraints/mathematical programs with variational equalities. For the value function of general parametric set-constrained programs, results on its Dini/directional derivative will be present in the forthcoming paper \cite{BY2022a}.}  

For any $\bar x,\bar y,u$, we denote by
{$$\mathcal C(\bar x,\bar y;u):=\left \{v\in \mathbb{L}(\bar x,\bar y;u) | f'_-((\bar x,\bar y);(u,v))\leq V_+'(\bar x;u),\  V_-'(\bar x;u)\leq f'_+((\bar x,\bar y);(u,v))\right \},$$}  the $y$-projection of the critical cone for the optimization problem
\begin{eqnarray*}
\min_{x,y} && f(x,y)-V(x) \nonumber\\
s.t. && P(x,y)\in \Gamma 
\end{eqnarray*} at $(x,y)=(\bar x,\bar y)$.
Given $u$ and $\alpha=1,0$,  we define
$$
		M^\alpha_u(\bar x,\bar y;\mathcal C(\bar x,\bar y;u)):=
		 \left \{( \zeta, \lambda) \left|\begin{array}{l}
			(\zeta,0) \in \alpha \partial f((\bar x,\bar y);(u,v))+\partial \langle P,\lambda\rangle ((\bar x,\bar y);(u,v)) , \\
			 \lambda \in \displaystyle \bigcup_{d\in DP(\bar x,\bar y)(u,v)}N_\Gamma(P(\bar x,\bar y); d), \quad v\in \mathcal C(\bar x,\bar y;u) \end{array} \right.\right \},
	$$ and 
	$$
		M^\alpha_0(\bar x,\bar y;\mathcal C(\bar x,\bar y;0)\cap \mathbb{S}):=
		 \left \{( \zeta, \lambda) \left|\begin{array}{l}
			(\zeta,0) \in \alpha \partial f((\bar x,\bar y);(0,v))+ \partial \langle P,\lambda\rangle ((\bar x,\bar y);(0,v)) , \\
			 \lambda \in \displaystyle \bigcup_{d\in DP(\bar x,\bar y)(0,v)}N_\Gamma(P(\bar x,\bar y); d), \quad  v\in \mathcal C(\bar x,\bar y;0)\cap \mathbb{S} \end{array} \right.\right \}.
	$$

Note that $M^\alpha_u(\bar x,\bar y;\mathcal C(\bar x,\bar y;u))$ and $M^\alpha_0(\bar x,\bar y;\mathcal C(\bar x,\bar y;0)\cap \mathbb{S})$ are   the set of  generalized Lagrange multipliers associated with the  problem 
$$\min_{x,y} f(x,y)\quad  \mbox{ s.t. }  P(x,y)\in \Gamma,  \quad {\bar x-x=0},$$ at $(\bar x, \bar y)$ in directions $\mathcal C(\bar x,\bar y;u)$ and $\mathcal C(\bar x,\bar y;0)\cap \mathbb{S}$, respectively.
 
We now present upper estimates for the limiting subdifferential of the value function.
\begin{thm}
\label{estimates} Let $ u\in \mathbb{R}^n$. \begin{itemize}
\item[{\rm (i)}]
Suppose that 
the restricted inf-compactness holds at $\bar x$ in direction $u$ with compact set $\Omega_u$, and the set-valued map $\Psi(x,y):=\Gamma-P(x,y)$ is  metrically subregular at $(\bar x, \tilde y, 0)$
for each $\tilde y\in S(\bar x;u)\cap\Omega_u$.
 {Moreover if $u\not =0$, assume that either $\Gamma$ is geometrically derivable at {$P(\bar x, \tilde y)$} for each $\tilde y\in S(\bar x;u)\cap\Omega_u$ or $P(x,y)$ is  directionally differentiable at $(\bar x, \tilde y)$ in direction $(u,v)$ for each $\tilde y\in S(\bar x;u)\cap\Omega_u$ and $v\in\mathbb L(\bar x,\tilde y;u)$.}
Then
	\begin{align} 
	\label{upperestim}
	\partial V(\bar x;u)\subseteq
	&\bigcup_{y\in  S(\bar x;u)\cap \Omega_u}
		 \left \{{\zeta}  \left  | (\zeta, \lambda) \in M^1_u(\bar x, y;  \mathcal C(\bar x,y;u))\cup M^1_0(\bar x, y; \mathcal C(\bar x,y;0)\cap\mathbb S)\right \} \right..
 	\end{align}
 	

\item[{\rm (ii)}] If in (i) the restricted inf-compactness at $\bar x$ in direction $u$ is replaced by the inner calmness* at $\bar x$ in direction $u$, 
then $M^1_0(\bar x, y; \mathcal C(\bar x,y;0)\cap\mathbb S)$ can be removed in (\ref{upperestim}) and hence  (\ref{upperestim}) becomes 
{\begin{align*}
	\nonumber
	\partial V(\bar x;u)\subseteq
	&\bigcup_{y\in  S(\bar x;u)}
	\left \{{\zeta}  \left  | (\zeta, \lambda) \in M^1_u(\bar x, y;  \mathcal C(\bar x,y;u))\right \} \right..
\end{align*}}

	\item[{\rm (iii)}] Suppose that there exists $\bar y\in S(\bar x)$ such that $S(x)$ is inner semicontinuous at $(\bar x,\bar y)$ in direction $u$ and the set-valued map $\Psi(x,y):=\Gamma-P(x,y)$ is  metrically subregular   at $(\bar x,\bar y,0)$. {Moreover if $u\not =0$, assume that  either $\Gamma$ is geometrically derivable at {$P(\bar x, \bar  y)$} or $P(x,y)$ is  directionally differentiable at $(\bar x, \bar  y)$ in direction $(u,v)$ for each  $v\in\mathbb L(\bar x,\bar y;u)$.} {Then the value function is continuous at $\bar x$ in direction $u$} and 
	\begin{align}
	\label{innerupere}
	\partial V(\bar x;u)\subseteq
	&	 \left \{\zeta \left | (\zeta, \lambda) \in M^1_u(\bar x, \bar y;  \mathcal C(\bar x,\bar y;u))\cup M^1_0(\bar x, \bar y; \mathcal C(\bar x,\bar y;0)\cap\mathbb S) \right. \right \}.
 	\end{align}
\item[{\rm (iv)}] {If in (iii), the inner semicontinuity of $S(x)$ at $(\bar x,\bar y)$ in direction $u$ is replaced by the inner calmness at $(\bar x,\bar y)$ in direction $u$, then $ M^1_0(\bar x, \bar y; \mathcal C(\bar x,\bar y;0)\cap\mathbb S)$ can be removed in (\ref{innerupere}) and 
hence (\ref{innerupere}) becomes
\begin{align*}
	\nonumber
	\partial V(\bar x;u)\subseteq
	&
	\left \{{\zeta}  \left  | (\zeta, \lambda) \in M^1_u(\bar x, \bar y;  \mathcal C(\bar x,\bar y;u))\right \} \right..
\end{align*}}
	\end{itemize}
\end{thm} 
\beginproof  
(i) 
Let $\zeta\in\partial V(\bar x;u)$. Then by Definition \ref{ads}, there exist sequences $t_k\downarrow0,\ u^k\rightarrow u,\ \zeta^k\rightarrow\zeta$ such that $V(\bar x+t_ku^k)\rightarrow V(\bar x)$
and 
$\zeta^k\in \widehat \partial V(\bar x+t_ku^k)$. It follows that 
$V(\bar x+t_ku^k)<V(\bar x)+\varepsilon$ for all $k$ large enough and hence by  the directional restricted inf-compactness, there exists $y^k\in S(\bar x+t_ku^k)\cap \Omega_u$.  Passing to a subsequence if necessary, we may assume that $y^k\rightarrow\tilde y$. Hence $\tilde y\in S(\bar x;u)\cap \Omega_u$.

For each $k$, since $\zeta^k\in \widehat \partial V(\bar x+t_ku^k)$, there exists a neighborhood ${\mathcal U}^k$ of $\bar x+t_ku^k$ satisfying
\begin{equation*}
V(x)-V(\bar x+t_ku^k)-\langle\zeta^k,x-(\bar x+t_ku^k)\rangle+\frac{1}{k}\|x-(\bar x+t_ku^k)\|\geq0\ \forall x\in{\mathcal U}^k.
\end{equation*}
It follows from the fact $V(x)=\displaystyle \inf_{y} \left \{f(x,y)+\delta_{ \Gamma}(P(x,y)) \right \}$ and $y^k\in S(\bar x+t_ku^k)$, that
\[
f(x,y)-\langle\zeta^k,x-(\bar x+t_ku^k)\rangle+\frac{1}{k}\|x-(\bar x+t_ku^k)\|+\delta_{ \Gamma}(P(x,y))\geq f(\bar x+t_ku^k,y^k),
\]
for any $(x,y)\in{\mathcal U}^k\times\mathbb R^m$.
Hence, the function 
$$
(x,y)\rightarrow f(x,y)-\langle\zeta^k,x-(\bar x+t_ku^k)\rangle+\frac{1}{k}\|x-(\bar x+t_ku^k)\|+\delta_{\Gamma}(P(x,y))$$ attains its local minimum at $(x,y)=(\bar x+t_ku^k,y^k)$. Then by the well-known Fermat's rule (see e.g., \cite[Proposition 1.30(i)]{Mor}) and the calculus rule in Proposition \ref{Calculus}, 
\begin{equation}
0\in \partial f(\bar x+t_ku^k,y^k)-(\zeta^k,0)+\frac{1}{k} {\mathcal B}\times \{0\}+\partial (\delta_{\Gamma}\circ P)(\bar x+t_ku^k,y^k).
\label{bounded}
\end{equation}
Now we consider two cases.\\
Case I: $\{\frac{ y^k-\tilde y}{t_k}\}$ is bounded. Passing to a subsequence if necessary, we can find $v\in \mathbb{R}^m$ such that $v^k:=\frac{ y^k-\tilde y}{t_k}\rightarrow v$ as $k\rightarrow \infty$. Since $y^k\in S(\bar x+t_ku^k)$ and $\tilde y\in S(\bar x)$, it follows that 
$$
 f'_-((\bar x,\tilde y);(u,v))\leq\lim_{k \rightarrow \infty}\frac{f(\bar x+t_ku^k,y^k)-f(\bar x,\tilde y)}{t_k}=\lim_{k\rightarrow \infty}\frac{V(\bar x+t_ku^k)-V(\bar x)}{t_k}\leq V_+'(\bar x;u).$$ Since $P(x,y)$ is Lipschitz continuous around $(\bar x,\tilde y)$, $P(\bar x+t_ku^k,y^k)\in \Gamma$ and $v^k\rightarrow v$,
 passing to a subsequence if necessary, we have
{
\[
\lim_{k \rightarrow \infty}\frac{P(\bar x+t_ku^k,\tilde y+t_k v^k)-P(\bar x,\tilde y)}{t_k}\in DP(\bar x,\tilde y)(u,v)\cap T_\Gamma(P(\bar x,\tilde y)) .
\]} 
Therefore
\begin{equation}\label{upperDini} f_-'((\bar x,\tilde y);(u,v)) \leq V_+'(\bar x;u), \qquad (u,v)\in \mathbb{L}(\bar x,\tilde y).
\end{equation}
But since 
$$ f(x,y)-V(x) \geq 0 \qquad \forall (x,y)\in {{\rm gph}{\mathcal F}}$$
and $\tilde y\in S(\bar x)$, $(\bar x,\tilde y)$ is an optimal solution to the above problem.
Hence by the first-order necessary optimality condition (\cite[Theorem 8.15]{RW}), {the lower Dini derivative of function $f(x,y)-V(x)$ at $(\bar x,\tilde y)$ along any feasible direction $(u,v)$ is nonnegative, i.e., $$(f-V)'_-((\bar x,\tilde y);(u,v))\geq0\quad \forall (u,v)\in T_{{\rm gph}\mathcal F}(\bar x,\tilde y).$$ Hence,} we have
$$f'_+((\bar x,\tilde y);(u,v)) +(- V)_+'(\bar x;u)\geq (f-V)'_-((\bar x,\tilde y);(u,v)) \geq 0 \qquad  \forall (u,v) \in T_{{\rm gph}{\mathcal F}}(\bar x, \tilde y). $$
{Since the metric subregularity  of $\Psi(x,y):=\Gamma-P(x,y)$ holds at $(\bar x,\tilde y,0)$, by Proposition \ref{Abadie}, under the assumptions of the theorem we have  
$T_{{\rm gph}{\mathcal F}}(\bar x, \tilde y)= \mathbb{L}(\bar x,\tilde y)$.}
Therefore {\begin{equation}\label{lowerDini} f'_+((\bar x,\tilde y);(u,v)) \geq  V_-'(\bar x;u), \qquad \forall v\in \mathbb{L}(\bar x,\tilde y;u).
\end{equation}} Combining (\ref{upperDini}) with (\ref{lowerDini}), we 
have $v\in \mathcal C(\bar x,\tilde y;u)$.
Taking limits as $k\rightarrow \infty$ in (\ref{bounded}), {by Proposition \ref{adsoc} since $y^k=\tilde y+t_k v^k$ and $v^k\rightarrow v$ we have} 
\begin{align}
0&\in\partial f((\bar x,\tilde y);(u,v))-(\zeta,0)+\partial (\delta_{\Gamma}\circ P)((\bar x,\tilde y);(u,v)).\label{eqn11}
\end{align} 

Since the metric subregularity of $\Psi(x,y)$ holds at $(\bar x,\tilde y,0)$ in direction $(u,v)$, by the directional chain rule in Proposition \ref{chainrule} we have
\begin{align}\partial (\delta_{\Gamma}\circ P)((\bar x,\tilde y);(u,v)) 
\subseteq \bigcup_{d\in DP(\bar x,\tilde y)(u,v)\cap T_\Gamma(P(\bar x,\tilde y))}\{\partial\langle P,\lambda\rangle((\bar x,\tilde y);(u,v))|\lambda\in N_{\Gamma}(P(\bar x,\tilde y); d)\}.
\label{eqn12}
\end{align} 
Combing (\ref{eqn11}) and (\ref{eqn12}), we obain
 $$(\zeta,0)\in\partial f((\bar x,\tilde y);(u,v))+\partial\langle P,\lambda\rangle((\bar x,\tilde y);(u,v)),$$ where  $\lambda\in N_{\Gamma}(P(\bar x,\tilde y); d)$ for some $d\in DP(\bar x,\tilde y)(u,v)\cap T_\Gamma(P(\bar x,\tilde y))$. This means that $(\zeta, \lambda) \in M^1_u(\bar x,\tilde y; \mathcal C(\bar x,\tilde y;u))$.

Case II: $\{\frac{ y^k-\tilde y}{t_k}\}$ is unbounded. Without loss of generality, assume $\lim_{k\rightarrow \infty} \frac{ \|y^k-\tilde y \|}{t_k}=\infty$.  Define $\tau_k:=\|y^k-\tilde y\|$. Then  $\frac{t_k}{\tau_k}\downarrow0$. Since the sequence $\{\frac{y^k-\tilde y}{\tau_k}\}$ is bounded, passing to a subsequence if necessary, we may assume that  there exist $v\in\mathbb S$ such that  $v^k:=\frac{ y^k-\tilde y}{\tau_k}\rightarrow v$. Define $\tilde u^k:=\frac{t_k}{\tau_k}u^k$. Then $\bar x+t_ku^k=\bar x+\tau_k\tilde u^k$ and $\tilde u^k\rightarrow0$.
  Since $y^k\in S(\bar x+\tau_k\tilde u^k)$ and $\tilde y\in S(\bar x)$, it follows that 
\begin{equation*}
f'_-((\bar x,\tilde y);(0,v))\leq\lim_{ k \rightarrow \infty}\frac{f(\bar x+\tau_k\tilde u^k,\tilde y+\tau_kv^k)-f(\bar x,\tilde y)}{\tau_k}=\lim_{ k \rightarrow \infty}\frac{V(\bar x+\tau_k\tilde u^k)-V(\bar x)}{\tau_k}\leq V'_+(\bar x;0).
\end{equation*}
Since $P(x,y)$ is Lipschitz continuous around $(\bar x,\tilde y)$, $P(\bar x+\tau_k\tilde u^k,y^k)\in \Gamma$, $\tilde u^k\rightarrow 0$ and $v^k\rightarrow v$,
 passing to a subsequence if necessary, we have
\[
\lim_{k \rightarrow \infty}\frac{P(\bar x+\tau_ku^k,\tilde y+\tau_k v^k)-P(\bar x,\tilde y)}{\tau_k}:=d\in DP(\bar x,\tilde y)(0,v)\cap T_\Gamma(P(\bar x,\tilde y)).
\]  Hence $v\in \mathbb{L}(\bar x,\tilde y;0).$
Therefore
\begin{equation*} f_-'((\bar x,\tilde y);(0,v)) \leq V_+'(\bar x;0), \qquad (0,v)\in \mathbb{L}(\bar x,\tilde y).
\end{equation*}
{Since 
$$ f(x,y)-V(x) \geq 0 \qquad \forall (x,y)\in {{\rm gph}{\mathcal F}}$$
and $\tilde y\in S(\bar x)$, $(\bar x,\tilde y)$ is an optimal solution to the above problem.
Hence by the first order necessary optimality condition, we have
$$f'_+((\bar x,\tilde y);(0,v)) +(- V)_+'(\bar x;0)\geq (f-V)'_-((\bar x,\tilde y);(0,v)) \geq 0 \qquad  \forall (0,v) \in T_{{\rm gph}{\mathcal F}}(\bar x, \tilde y). $$
Since the metric subregularity  of $\Psi(x,y):=\Gamma-P(x,y)$ holds at $(\bar x,\tilde y,0)$, by Proposition \ref{Abadie}, under the assumptions of the theorem we have  
	$T_{{\rm gph}{\mathcal F}}(\bar x, \tilde y)= \mathbb{L}(\bar x,\tilde y)$.
Therefore \begin{equation*} f'_+((\bar x,\tilde y);(0,v)) \geq  V_-'(\bar x;0), \qquad \forall v\in \mathbb{L}(\bar x,\tilde y;0).
\end{equation*}
In summary, we have
$$f'_-((\bar x,\tilde y);(0,v))\leq V'_+(\bar x;0),\ V_-'(\bar x;0)\leq f'_+((\bar x,\tilde y);(0,v))\qquad \forall v\in \mathbb{L}(\bar x,\tilde y;0).$$  Hence $v\in \mathcal C(\bar x,\tilde y;0)$.}
Taking limits as $k\rightarrow \infty$ in (\ref{bounded}), since $\bar x+t_ku^k=\bar x+\tau_k\tilde u^k$ and $\tilde u^k\rightarrow0$, following a similar process as in Case I we have 
\begin{align*}
0&\in\partial f((\bar x,\tilde y);(0,v))-(\zeta,0)+\partial (\delta_{\Gamma}\circ P)((\bar x,\tilde y);(0,v)).
\end{align*} 
Hence there exists $\lambda \in N_{\Gamma}(P(\bar x,\tilde y); d)$ with $d\in DP(\bar x,\tilde y)(0,v)\cap T_\Gamma(P(\bar x,\tilde y))$ and $$(\zeta,0)\in \partial f((\bar x,\tilde y);(0,v))+\partial\langle P,\lambda \rangle((\bar x,\tilde y);(0,v)).$$ This means that $(\zeta, \lambda) \in M^1_0(\bar x,\tilde y; \mathcal C(\bar x,\tilde y;0)\cap \mathbb{S})$. 

The proof of (i) is complete by combing cases I and II.

{
(ii) When $S(x)$ is inner calm* at $\bar x$ in direction $u$, passing to a subsequence if necessary, one can choose in the above proof $y^k\in S(\bar x+t_ku^k)$ such that $\{(y^k-\bar y)/t_k\}$ is bounded. Consequently, Case II is impossible, and {$M^1_0(\bar x, y; \mathcal C(\bar x,y;0)\cap\mathbb S)$ is not needed.}
}

(iii) When $S(x)$ is inner semicontinuous at some point $\bar y\in S(\bar x)$ in direction $u$, the restricted inf-compactness holds at $\bar x$ in direction $u$ and $\bar y\in S(\bar x;u)$. One can prove the result by taking $\tilde y=\bar y$ in the proof of (i). {And by Proposition \ref{directional-conti}, $V(x)$ is continuous at $\bar x$ in direction $u$.}

{
(iv) When $S(x)$ is inner calm at $(\bar x,\bar y)$ in direction $u$, passing to a subsequence if necessary, one can choose in the above proof $y^k\in S(\bar x+t_ku^k)$ satisfying that $y^k\rightarrow\bar y$ and $\{(y^k-\bar y)/t_k\}$ is bounded. Consequently, Case II is impossible, and {$M^1_0(\bar x, \bar y; \mathcal C(\bar x,y;0)\cap\mathbb S)$ is not needed.}
}
\endproof

\begin{remark}
\label{remark3.1}

{When $u=0$, $$M_0^\alpha(x,y; \mathcal C(x,y;0)\cap\mathbb S)\subseteq M_0^\alpha(x,y; \mathcal C(x,y;0))=M_u^\alpha(x,y; \mathcal C(x,y;u)).$$ Since $0\in \mathcal C(x,y;0)$ {and $0\in DP(x,y)(0,0)$} always, we have
{$$\bigcup_{d\in DP(x,y)(0,v),\atop v\in  \mathcal C(x,y;0)}N_\Gamma(P(x,y); d) =N_\Gamma(P(x,y)),$$}
and therefore we have
\begin{eqnarray*}
\lefteqn{M_u^\alpha(x,y; \mathcal C(x,y;0))\cup M_0^\alpha(x,y; \mathcal C(x,y;0)\cap\mathbb S)=M_0^\alpha(x,y; \mathcal C(x,y;0))=M^\alpha(x,y)}\\
& &:=\left \{(\zeta, \lambda) |(\zeta,0)\in \alpha \partial f(x,y)+ \partial \langle P,\lambda \rangle (x, y), \lambda\in N_\Gamma\left (P(x,y)\right ) \right \}.\end{eqnarray*}
Hence even when $u=0$, Theorem \ref{estimates}(i) improves \cite[Theorem 3.6]{LY,LYErru} 
in that the metric subregularity is assumed which are weaker than  NNAMCQ as required in \cite[Theorem 3.6]{LY,LYErru}.}
\end{remark}
{Under the assumptions in Theorem \ref{estimates}, the value function may not be Lipschitz continuous, even with respect to a directional neighborhood. This means that $\partial^\infty V(x;u)$ may contain nonzero elements and  so it is also meaningful to give an estimate for $\partial^\infty V(x;u)$.} Based on the analysis in Remark \ref{remark3.1}, even when $u=0$, Theorem \ref{isestimates}(i) improves \cite[Theorem 3.6]{LY,LYErru} 
in that the metric subregularity is assumed which is weaker than  NNAMCQ as required in \cite[Theorem 3.6]{LY,LYErru}.
\begin{thm}
	\label{isestimates}  Let $u\in \mathbb{R}^n$.
	\begin{itemize}
	\item[{\rm (i)}]
	Suppose that  the restricted inf-compactness holds at $\bar x$ in direction $u$ with compact set $\Omega_u$.  Furthermore suppose that  the set-valued map $\Psi(x,y):=\Gamma-P(x,y)$ is  metrically subregular at $(\bar x, y)$ for each $y\in S(\bar x;u)\cap\Omega_u$. {Moreover if $u\not =0$, suppose either $\Gamma$ is geometrically derivable
	at {$P(\bar x, \tilde y)$} for each $\tilde y\in S(\bar x;u)\cap\Omega_u$ or $P(x,y)$ is  directionally differentiable at $(\bar x,\tilde y)$  in direction $(u,v)$ for each $\tilde y\in S(\bar x;u)\cap\Omega_u$ and $v\in\mathbb L(\bar x,\tilde y;u)$.}
	Then
	\begin{align}\label{upperestimaym}
	\partial^\infty V(\bar x;u)\subseteq
	&\bigcup_{y\in  S(\bar x;u)\cap \Omega_u}	 \left \{\zeta : (\zeta, \lambda) \in M_u^0(\bar x,y; \mathcal C(\bar x,y;u))\cup M_0^0(\bar x,y; \mathcal C(\bar x,y;0)\cap\mathbb S) \right \}.
	\end{align}

\item[{\rm (ii)}] 
{If in (i) the restricted inf-compactness at $\bar x$ in direction $u$ is replaced by the inner calmness* at $\bar x$ in direction $u$, then $M^0_0(\bar x, y; \mathcal C(\bar x,y;0)\cap\mathbb S)$ can be removed in (\ref{upperestimaym}) and hence (\ref{upperestimaym}) becomes} 
{\begin{align*}
	\nonumber
	\partial^\infty V(\bar x;u)\subseteq
	&\bigcup_{y\in  S(\bar x;u)}
	\left \{{\zeta}  \left  | (\zeta, \lambda) \in M^0_u(\bar x, y;  \mathcal C(\bar x,y;u))\right \} \right..
\end{align*}}

	\item[{\rm (iii)}] Suppose that there exists $\bar y\in S(\bar x)$ such that $S(x)$ is inner semicontinuous at $(\bar x,\bar y)$ in direction $u$ and the set-valued map $\Psi(x,y):=\Gamma-P(x,y)$ is  metrically subregular   at $(\bar x,\bar y,0)$. {Moreover if $u\not =0$, assume that  either $\Gamma$ is geometrically derivable at {$P(\bar x, \bar  y)$} or $P(x,y)$ is  directionally differentiable at $(\bar x, \bar  y)$ in direction $(u,v)$ for each  $v\in\mathbb L(\bar x,\bar y;u)$.}  Then 
	\begin{align}
	\label{innerupereaym}
	\partial^\infty V(\bar x;u)\subseteq
	& \left \{\zeta : ( \zeta,\lambda) \in M_u^0(\bar x,\bar y; \mathcal C(\bar x,\bar y;u))\cup M_0^0(\bar x,\bar y; \mathcal C(\bar x,\bar y;0)\cap\mathbb S) \right \}.
	\end{align}

\item[{\rm (iv)}]{If in (iii), the inner semicontinuity of $S(x)$ at $(\bar x,\bar y)$ in direction $u$ is replaced by the inner calmness at $(\bar x,\bar y)$ in direction $u$, then $ M^0_0(\bar x, \bar y; \mathcal C(\bar x,\bar y;0)\cap\mathbb S)$ can be removed in (\ref{innerupereaym}) and hence (\ref{innerupereaym}) becomes}
\begin{align*}
	\nonumber
	\partial^\infty V(\bar x;u)\subseteq
	&
	\left \{{\zeta}  \left  | (\zeta, \lambda) \in M^0_u(\bar x, \bar y;  \mathcal C(\bar x,\bar y;u))\right \} \right..
\end{align*}
	\end{itemize}
\end{thm} 
\beginproof (i)
Let $\zeta\in\partial^\infty V(\bar x;u)$. Then by Definition \ref{ads}, there exist sequences $l_k\downarrow0,\ t_k\downarrow0,\ u^k\rightarrow u,\ \zeta^k\rightarrow\zeta$ such that $V(\bar x+t_ku^k)\rightarrow V(\bar x)$
and 
$\zeta^k\in \widehat \partial V(\bar x+t_ku^k)$ with $l_k\zeta^k\rightarrow\zeta$. 

%
Following a similar process as in the proof of Theorem \ref{estimates}, we can obtain $y^k \in S(\bar x+ t_k u^k) \cap \Omega_u$ and $\tilde y \in S(\bar x;u)\cap \Omega_u$ such that $y^k\rightarrow \tilde y$. Moreover (\ref{bounded}) holds. Multiplying both sides of (\ref{bounded}) by $l_k$, we have
\begin{equation}
0\in l_k\partial f(\bar x+t_ku^k,y^k)-l_k(\zeta^k,0)+\frac{l_k}{k}{\mathcal B}\times \{0\}+\partial (\delta_{\Gamma}\circ P)(\bar x+t_ku^k,y^k).\label{bounded2}
\end{equation}
Now we consider two cases.\\
Case I: $\{\frac{ y^k-\tilde y}{t_k}\}$ is bounded. Passing to a subsequence if necessary, there exists $v\in \mathbb{R}^m$ such that $\frac{ y^k-\tilde y}{t_k}\rightarrow v$. 
Similarly as in the proof of Theorem \ref{estimates}(i), we have
$$f'_-((\bar x, \tilde y);(u,v))\leq V_+'(\bar x;u),\ V_-'(\bar x;u)\leq f'_+((\bar x, \tilde y);(u,v)) \mbox{ and } DP(\bar x,\tilde y)(u,v) \cap T_C(P(\bar x,\tilde y))\neq\emptyset.$$

Taking limits as $k\rightarrow \infty$ in (\ref{bounded2}), {by Proposition \ref{adsoc} we have} 
\begin{align*}
0&\in-(\zeta,0)+\partial (\delta_{\Gamma}\circ P)((\bar x,\tilde y);(u,v)).
\end{align*} 

 By assumption, the metric subregularity of $\Psi(x,y):=\Gamma-P(x,y)$ holds at $(\bar x,\tilde y,0)$ in direction $(u,v)$. Hence by Proposition \ref{chainrule} we have
\[\partial (\delta_{\Gamma}\circ P)((\bar x,\tilde y);(u,v)) 
\subseteq \bigcup_{d\in DP(\bar x,\tilde y)(u,v)\cap T_\Gamma(P(\bar x,\tilde y))}\left \{\partial\langle P,\lambda \rangle((\bar x,\tilde y);(u,v))|\lambda\in N_{\Gamma}(P(\bar x,\tilde y); d)\right \}.
\]
Hence $$(\zeta,0) \in \partial\langle\lambda,P\rangle((\bar x,\tilde y);(u,v)), \ \lambda\in N_\Gamma(P(\bar x,\tilde y);d),\ d\in DP(\bar x,\tilde y)(u,v)\cap T_\Gamma(P(\bar x,\tilde y)).$$ This means that $(\zeta,\lambda)\in M^0_u(\bar x,\tilde y;\mathcal C(\bar x,\tilde y;u))$.

Case II: $\{\frac{ y^k-\tilde y}{t_k}\}$ is unbounded. Without loss of generality, assume $\lim_{k\rightarrow \infty} \frac{ \|y^k-\tilde y \|}{t_k}=\infty$.  Define $\tau_k:=\|y^k-\tilde y\|$. Then  $\frac{t_k}{\tau_k}\downarrow0$. Since the sequence $\{\frac{y^k-\tilde y}{\tau_k}\}$ is bounded, without loss of generality, assume there exist $v\in\mathbb S$ and a sequence $v^k\rightarrow v$ such that $ y^k=\tilde y+\tau_kv^k$. Define $\tilde u^k:=\frac{t_k}{\tau_k}u^k$. Then $\bar x+t_ku^k=\bar x+\tau_k\tilde u^k$ and $\tilde u^k\rightarrow0$. Similarly as in the proof of Case II of Theorem \ref{estimates}(i), we have
$$ f'_-((\bar x, \tilde y);(0,v))\leq V'_+(\bar x;0), V'_-(\bar x;0)\leq f'_+((\bar x, \tilde y);(0,v))\mbox{ and } DP(\bar x,\tilde y)(0,v) \cap T_\Gamma(P(\bar x,\tilde y))\neq\emptyset.$$
Taking limits as $k\rightarrow \infty$ in (\ref{bounded2}), following a similar process as in Case I we have 
\begin{align*}
0&\in-(\zeta,0)+\partial (\delta_{\Gamma}\circ P)(\bar x,\tilde y;0,v).
\end{align*} 
Hence there exists $\lambda\in N_{\Gamma}(P(\bar x,\tilde y);d)$ with $d\in DP(\bar x,\tilde y)(0,v) \cap T_\Gamma(P(\bar x,\tilde y))$ and $$(\zeta,0) \in \partial\langle\lambda,P\rangle((\bar x,\tilde y);(0,v)), \ \lambda\in N_\Gamma(P(\bar x,\tilde y);DP(\bar x,\tilde y)(0,v)).$$ 
Hence $(\zeta,\lambda)\in M_0^0(\bar x,\tilde y;\mathcal C(\bar x,\tilde y;0)\cap \mathbb{S})$.

{
	(ii) {When $S(x)$ is inner calm* at $\bar x$ in direction $u$,  passing to a subsequence if necessary, one can choose in the above proof $y^k\in S(\bar x+t_ku^k)$ satisfying that $\{(y^k-\bar y)/t_k\}$ is bounded. Consequently, Case II is impossible, and $M^0_0(\bar x, y; \mathcal C(\bar x,y;0)\cap\mathbb S)$ is not needed.}
}

(iii)   When $S(x)$ is inner semicontinuous at some point $\bar y\in S(\bar x)$ in direction $u$, the restricted inf-compactness holds at $\bar x$ in direction $u$ and $\bar y\in S(\bar x;u)$. In the proof of (i), taking $\tilde y=\bar y$, one obtains the desired result. 

{
	(iv) When $S(x)$ is inner calm at $(\bar x,\bar y)$ in direction $u$, the restricted inf-compactness holds at $\bar x$ in direction $u$. Passing to a subsequence if necessary, one can choose in the above proof $y^k\in S(\bar x+t_ku^k)$ satisfying that $y^k\rightarrow\bar y$ and $\{(y^k-\bar y)/t_k\}$ is bounded. Consequently, Case II is impossible, and 
	{$M^0_0(\bar x, y; \mathcal C(\bar x,\bar y;0)\cap\mathbb S)$ is not needed.}
}
\endproof

From items (ii), (iv) in Theorems \ref{estimates} and \ref{isestimates}, when the optimal solution map $S(x)$ satisfies the directional inner calmness*/calmness, the set $M^\alpha_0(\bar x,y;\mathcal C(\bar x,y;0)\cap\mathbb S)(\alpha=0,1)$ can be removed from the estimates.  Usually when $S(x)$ does not satisfy the directional inner calmness/calmness* conditions, $M^\alpha_0(\bar x,y;\mathcal C(\bar x,y;0)\cap\mathbb S)(\alpha=0,1)$ is needed to estimate the subdifferential of the value function. 
In the following example, the value function $V(x)$ is not lipschitz, $M^0_0(\bar x,y;\mathcal C(\bar x,y;0)\cap\mathbb S)$ has a nonzero element and $M^0_u(\bar x,\bar y;\mathcal C(\bar x,\bar y;u))=\emptyset$  for $u\not =0$. Hence $M^0_0(\bar x,y;\mathcal C(\bar x,y;0)\cap\mathbb S)$ is needed to show that the value function is not lipschitz continuous at $\bar x$ in direction $u$.

\begin{example}
	\begin{align*}
		\min_y\ y\ \  {\rm s.t.}\ (x,y) \in \Gamma:=\{(x,y)|x-y^3\leq0\}.
	\end{align*}	
\end{example}

In this example, $f(x,y):=y, P(x,y):=(x,y)$. One can easily find that $V(x)=\sqrt[3]{x}$ and $S(x)=\sqrt[3]{x}$.  Consider the point $(\bar x,\bar y)=(0,0)$ and direction $u=1$. Then $S(x)$ as a single-valued continuous function is inner semicontinuous at $(\bar x,\bar y)$ in direction $u$. Obviously,  $V(x)$ is not Lipschitz continuous at $\bar x$ in direction $u$, and the inner calmness at $(\bar x,\bar y)$/the inner calmness* at $\bar x$ in direction $u$ does not hold for $S(x)$. $\partial V(\bar x;u)=\emptyset$ and $\partial^\infty V(\bar x;u)=\mathbb{R}_+$. By calculation, 
\begin{eqnarray*}
T_\Gamma(0,0)&=& \{0\}\times \mathbb{R},\\
N_\Gamma(0,0; (0,v))&=& \mathbb{R}_+\times \{0\}, \ \mbox{ for } v\in \mathbb{R},\\
\mathbb{L}(0,0;0)&=& \{v|(u,v)\in T_\Gamma(0,0)\}=\mathbb{R},\\
\mathbb{L}(0,0;u)&=& \emptyset  \ \mbox{ for } u\not =0,\\
\mathcal C(0,0;0)&=& {\{ v| -\infty \leq v \leq\infty\}=\mathbb{R}},\\
\mathcal C(0,0;u)&=&\emptyset  \ \mbox{ for } u\not =0,\\
M^\alpha_0(0,0;\mathcal C(0,0;0)\cap \mathbb S)&=& \left \{(\zeta,\lambda)\mid \begin{array}{l} (\zeta,0)= (0, \alpha v) +\lambda,\\
 \lambda \in N_\Gamma(0,0; (0,v)), v\in \mathcal C(0,0;0)\cap \mathbb S \end{array} \right \},\\
M^\alpha_u(0,0;\mathcal C(0,0;u))&=& \left \{(\zeta,\lambda)\mid \begin{array}{l} (\zeta,0)= (0, \alpha v) +\lambda,\\ \lambda \in N_\Gamma(0,0; (0,v)), v\in \mathcal C(0,0;u)\end{array} \right \}, \mbox{ for } u\not =0.
\end{eqnarray*}
 Hence, for $u\not =0$, $M^1_u(\bar x,\bar y;\mathcal C(\bar x,\bar y;u))=M^0_u(\bar x,\bar y;\mathcal C(\bar x,\bar y;u))=M^1_0(\bar x,\bar y;\mathcal C(\bar x,\bar y;0)\cap\mathbb S)=\emptyset$ and $M^0_0(\bar x,\bar y;\mathcal C(\bar x,\bar y;0)\cap\mathbb S)=\mathbb{R}_+$. Then the conclusions in Theorems 3.1 and 3.2 hold since one has $\partial V(\bar x;u)=\emptyset$ and $\partial^\infty V(\bar x;u)\subseteq \mathbb{R}_+$. 


\begin{remark}\label{remark3.2} Note that Long et al. \cite[Theorem  5.11]{Long} obtained some upper estimates for the directional limiting and singular subdifferential of a constrained optimization problem  under a stronger version of the directional inner semicontinuity (\cite[Definition 4.4(i)]{Long} of $S(x)$. Since our assumptions are weaker than that in \cite[Theorems 5.10 and 5.11]{Long}, our results can not be obtained by using \cite[Theorem 5.11]{Long}.

Recall that $\mathcal F(x):=\{y|P(x,y)\in\Gamma\}$. Then it is obvious that the value function can be rewritten as
$$V(x)=\inf_y \vartheta (x,y), \mbox{ where }
\vartheta (x,y):=f(x,y)+\delta_{{\rm gph}{\mathcal F}}(x,y).$$
Long et al. \cite[Theorem  5.10]{Long} obtained some upper estimates for the directional limiting and singular subdifferential of the value function $V(x)$ in terms of the corresponding  directional limiting and singular subdifferential for $\vartheta(x,y)$ under a stronger version of the directional inner semicontinuity. Benko et al. \cite{BHA} obtained upper estimates for the  directional limiting subdifferential $\tilde\partial V(x;(u,\xi))$ for $V(x)$ in directions $u$ and $\xi$ from $\bar x$ and $V(\bar x)$, respectively. Using the relationship between the subdifferentials, $\tilde\partial V(x;(u,\xi))$ and $\partial V(x;u)$ obtained in \cite[Corollary 4.1]{BHA}, one can derive some upper estimate for $\partial V(x;u)$. However Theorem \ref{estimates} can not be derived in this way from \cite{BHA}. First, observing that in the proof of Theorem \ref{estimates}, to state the upper estimates in terms of problem data $P(x,y)$ and $\Gamma$, calculus rules like Proposition \ref{Abadie} are needed. Secondly, compared with \cite[Theorem 4.2]{BHA}, Theorem \ref{estimates} assumes the  directional restricted inf-compactness which is weaker than its non-directional counterpart, and  replaces the solution set $S(\bar x)$ by its subset $S(\bar x;u)$ which makes the upper estimates sharper.
\end{remark}

Combining Propositions \ref{VLip} and \ref{Vlip2}   and Theorems \ref{estimates} and  \ref{isestimates}, we obtain the following sufficient conditions for the directional Lipschitz continuity of value function $V(x)$.
\begin{thm}\label{Prop3.4} Let $u\in \mathbb{R}^n$.
	\begin{itemize}
	\item[{\rm (i)}]
	Suppose that  the restricted inf-compactness holds at $\bar x$ in direction $u$ with compact set $\Omega_u$. Furthermore suppose $V(x)$ is continuous at $\bar x$ in direction $u$ if $u\not =0$ and    the set-valued map $\Psi(x,y):=\Gamma-P(x,y)$ is  metrically subregular at $(\bar x, \tilde y)$ for each $\tilde y\in S(\bar x;u)\cap\Omega_u$.  {Moreover if $u\not =0$, suppose either $\Gamma$ is geometrically derivable at {$P(\bar x,\tilde y)$} for each $\tilde y\in S(\bar x;u)\cap\Omega_u$ or $P(x,y)$ is directionally differentiable in direction $(u,v)$ for each $\tilde y\in S(\bar x;u)\cap\Omega_u$ and $v\in\mathbb L(\bar x,\tilde y;u)$.}
	If 
	\begin{align*}
	\bigcup_{y\in  S(\bar x;u)\cap \Omega_u}&  \left \{\zeta \left | (\zeta,\lambda) \in M_u^0(\bar x,y; \mathcal C(\bar x,y;u))\cup M_0^0(\bar x,y; \mathcal C(\bar x,y;0)\cap\mathbb S) \right .\right \}
	=\{0\},
	\end{align*}
	then $V(x)$ is Lipschitz around $\bar x$ in direction $u$, and
	\begin{align*}
	\nonumber
	\emptyset \not =\partial V(\bar x;u)\subseteq
	&	\bigcup_{y\in  S(\bar x;u)\cap \Omega_u}  \left \{\zeta \left | (\zeta, \lambda) \in M_u^1(\bar x,y; \mathcal C(\bar x,y;u))\cup M_0^1(\bar x,y; \mathcal C(\bar x,y;0)\cap\mathbb S) \right  .\right \}.
\end{align*}	

\item[{\rm (ii)}]
{If in (i) the restricted inf-compactness at $\bar x$ in direction $u$ is replaced by the inner calmness* at $\bar x$ in direction $u$, then $M^\alpha_0(\bar x, y; \mathcal C(\bar x,y;0)\cap\mathbb S)(\alpha=0,1)$ can be removed in conditions of {\rm (i)} and one has, } 
if 
{\begin{align*}
	\bigcup_{y\in  S(\bar x;u)}&  \left \{\zeta \left | (\zeta,\lambda) \in M_u^0(\bar x,y; \mathcal C(\bar x,y;u)) \right .\right \}
	=\{0\},
\end{align*}
then $V(x)$ is Lipschitz around $\bar x$ in direction $u$, and
\begin{align*}
	\nonumber
	\emptyset \not =\partial V(\bar x;u)\subseteq
	&	\bigcup_{y\in  S(\bar x;u)}  \left \{\zeta \left | (\zeta, \lambda) \in M_u^1(\bar x,y; \mathcal C(\bar x,y;u)) \right  .\right \}.
\end{align*}}

	\item[{\rm (iii)}] Suppose that  there exists $\bar y\in S(\bar x)$ such that $S(x)$ is inner semicontinuous at $(\bar x,\bar y)$ in direction $u$  and  the set-valued map $\Psi(x,y):=\Gamma-P(x,y)$ is  metrically subregular   at $(\bar x,\bar y,0)$.
	Moreover if $u\not =0$, suppose either $\Gamma$ is geometrically derivable at $P(\bar x,\bar y)$ or $P(x,y)$ is  directionally differentiable at $(\bar x,\bar y)$ in direction $(u,v)$ for each $v\in\mathbb L(\bar x,\bar  y;u)$. If 
	\begin{align*}
	&	 \left \{\zeta \left | (\zeta, \lambda) \in M_u^0(\bar x,\bar y; \mathcal C(\bar x,\bar y;u))\cup M_0^0(\bar x,\bar y; \mathcal C(\bar x,\bar y;0)\cap\mathbb S) \right .\right \}=\{0\},
	\end{align*}
	then $V(x)$ is Lipschitz around $\bar x$ in direction $u$, and 
	\begin{align*}
	\nonumber
	\emptyset \not =\partial V(\bar x;u)\subseteq
	&	 \left \{\zeta \left | (\zeta, \lambda) \in M_u^1(\bar x,\bar y; \mathcal C(\bar x,\bar y;u))\cup M_0^1(\bar x,\bar y; \mathcal C(\bar x,\bar y;0)\cap\mathbb S) \right. \right \}.
 	\end{align*}

	\item[{\rm (iv)}] 
{If in (iii), the inner semicontinuity of $S(x)$ at $(\bar x,\bar y)$ in direction $u$ is replaced by the inner calmness at $(\bar x,\bar y)$ in direction $u$, then $ M^\alpha_0(\bar x, \bar y; \mathcal C(\bar x,\bar y;0)\cap\mathbb S)(\alpha=0,1)$ can be removed in conditions of {\rm (iii)} and one has, }
if 
\begin{align*}
	&	 \left \{\zeta \left | (\zeta, \lambda) \in M_u^0(\bar x,\bar y; \mathcal C(\bar x,\bar y;u)) \right .\right \}=\{0\},
\end{align*}
then $V(x)$ is Lipschitz around $\bar x$ in direction $u$, and 
\begin{align*}
	\nonumber
	\emptyset \not =\partial V(\bar x;u)\subseteq
	&	 \left \{\zeta \left | (\zeta, \lambda) \in M_u^1(\bar x,\bar y; \mathcal C(\bar x,\bar y;u)) \right. \right \}.
\end{align*}

	\end{itemize}
\end{thm}
\beginproof  (i) Since  the restricted inf-compactness holds at $\bar x$ in direction $u$ and the value function is continuous at $\bar x$ in direction $u$, 
by Proposition \ref{VLip},  $V(x)$ is Lipschitz around $\bar x$ in direction $u$ if and only if $\partial^\infty V(\bar x; u)=\{0\}$. By the upper estimate for $\partial^\infty V(\bar x; u)$ in Theorem \ref{isestimates}(i), we have $\partial^\infty V(\bar x; u)=\{0\}$.
Hence $V(x)$ is Lipschitz at $\bar x$ in direction $u$.

{
	(ii) When $S(x)$ is inner calm* at $\bar x$ in direction $u$ and the value function is continuous at $\bar x$ in direction $u$ if $u\not =0$, { by Proposition \ref{Vlip2},  $V(x)$ is Lipschitz around $\bar x$ in direction $u$ if and only if $\partial^\infty V(\bar x; u)=\{0\}$.  By the upper estimate for $\partial^\infty V(\bar x; u)$ in Theorem \ref{isestimates}(ii), we have $\partial^\infty V(\bar x; u)=\{0\}$.}
	Hence $V(x)$ is Lipschitz at $\bar x$ in direction $u$.
}

(iii) Since $S(x)$ is inner semicontinuous  $(\bar x, \bar y) $ in direction $u$, the restricted inf-compactness holds at $\bar x$ in direction $u$ and the value function is continuous at $\bar x$ in direction $u$. Hence by Proposition \ref{VLip},  $V(x)$ is Lipschitz around $\bar x$ in direction $u$ if and only if $\partial^\infty V(\bar x; u)=\{0\}$. By the upper estimate for $\partial^\infty V(\bar x; u)$ in Theorem \ref{isestimates}(ii), we have $\partial^\infty V(\bar x; u)=\{0\}$. Hence  $V(x)$ is Lipschitz around $\bar x$ in direction $u$.

{
(iv) Since $S(x)$ is inner calm  $(\bar x, \bar y) $ in direction $u$, the restricted inf-compactness holds at $\bar x$ in direction $u$ and the value function is continuous at $\bar x$ in direction $u$. Hence by Proposition \ref{VLip},  $V(x)$ is Lipschitz around $\bar x$ in direction $u$ if and only if $\partial^\infty V(\bar x; u)=\{0\}$. By the upper estimate for $\partial^\infty V(\bar x; u)$ in Theorem \ref{isestimates}(iv), we have $\partial^\infty V(\bar x; u)=\{0\}$. Hence  $V(x)$ is Lipschitz around $\bar x$ in direction $u$.
}
\endproof

\section{Application to special cases}
In this section we apply our results to various important special cases. For these special cases, some of the assumptions are automatically satisfied and the expressions are simpler.
\subsection{Smooth case}
In this section, we consider the case where  all functions $f,P$ are smooth. In this case,
{\begin{eqnarray*}
\mathbb{L}(x,y;u)&=&\left \{v\left | \nabla P(x,y)(u,v)\in T_\Gamma(P(x,y))
	\right .\right \},\\
\mathcal C(x,y;u)&=& \left \{v\in \mathbb{L}(x,y;u) |  V_-'(x;u)\leq   \nabla f(x,y)(u,v)\leq V_+'(x;u)\right \}. 
\end{eqnarray*}} 
{For $\alpha=0,1$, we define the set of generalized Lagrange multipliers associated with $y\in{S}(x)$  as
\begin{eqnarray*}
\lefteqn{\Lambda^\alpha_u(x,y; \mathcal C(x,y;u))  :=}\\
&& \left \{\lambda\left|0= \alpha\nabla_y f(x,y)+ \nabla_y P(x,y)^T \lambda, \  \lambda\in \bigcup_{ v\in \mathcal C(x,y;u)} N_\Gamma\left (P(x,y); \nabla P(x,y)(u,v) \right )\right. \right \}, 
\end{eqnarray*}
and
\begin{eqnarray*}
\lefteqn{\Lambda^\alpha_0(x,y; \mathcal C(x,y;0)\cap \mathbb{S})  :=}\\
&& \left \{\lambda\left|0= \alpha\nabla_y f(x,y)+ \nabla_y P(x,y)^T \lambda, \  \lambda\in \bigcup_{ v\in \mathcal C(x,y;0)\cap \mathbb{S}} N_\Gamma\left (P(x,y); \nabla P(x,y)(0,v) \right )\right. \right \}, 
\end{eqnarray*} in directions $\mathcal C(x,y;u)$ and $\mathcal C(x,y;0)\cap \mathbb{S}$
respectively. Then we have
{\begin{eqnarray*}
\lefteqn{M^\alpha_u( x, y;\mathcal C( x, y;u)) := }\\
 && \left \{(\zeta,\lambda ) \left|\zeta  = \alpha\nabla_x f( x,y)+ \nabla_x P( x,y)^T \lambda, \  \lambda\in \Lambda^\alpha_u( x,y; \mathcal C( x, y;u))\right. \right \}\\
\lefteqn{M^\alpha_0( x, y;\mathcal C( x, y;0)\cap \mathbb{S}) := }\\
&& \left \{(\zeta,\lambda ) \left|\zeta= \alpha\nabla_x f( x,y)+ \nabla_x P( x,y)^T \lambda, \  \lambda\in \Lambda^\alpha_0( x,y; \mathcal C( x, y;0)\cap\mathbb S)\right. \right \}.
\end{eqnarray*}}}
  Hence Theorem \ref{Prop3.4} has the following consequence. 
\begin{prop}\label{smooth} Assume that all functions $f,P$ are smooth. Let $u\in \mathbb{R}^n$.
	\begin{itemize}
	\item[{\rm (i)}]
	Suppose that  the restricted inf-compactness holds at $\bar x$ in direction $u$ with compact set $\Omega_u$. Furthermore suppose that the value function $V(x)$ is continuous at $\bar x$ in direction $u$ if $u\not =0$ and   the set-valued map $\Psi(x,y):=\Gamma-P(x,y)$ is  metrically subregular at $(\bar x, \tilde y,0)$ for all $\tilde y\in S(\bar x;u)\cap \Omega_u$.  
	If 
	\begin{align*}
	\bigcup_{y\in  S(\bar x;u)\cap \Omega_u}& 	\left \{\nabla_x P(\bar x, y)^T \lambda \left | \lambda \in \Lambda_u^0(\bar x,y; \mathcal C(\bar x,y;u))\cup\Lambda_0^0(\bar x,y; \mathcal C(\bar x,y;0)\cap\mathbb S) \right .\right \}
	=\{0\},
	\end{align*}
	then $V(x)$ is Lipschitz around $\bar x$ in direction $u$, and
	\begin{align*}
	\nonumber
	\emptyset &\not =\partial V(\bar x;u)\\
	&\subseteq
		\bigcup_{y\in  S(\bar x;u)\cap \Omega_u}  \left \{\nabla_x f(\bar x, y)+\nabla_x P(\bar x,  y)^T\lambda \left |\lambda \in \Lambda_u^1(\bar x,y; \mathcal C(\bar x,y;u))\cup \Lambda_0^1(\bar x,y; \mathcal C(\bar x,y;0)\cap\mathbb S)\right . \right \}.
 	\end{align*}

\item[{\rm (ii)}]
{Suppose that  the inner calmness*  holds at $\bar x$ in direction $u$. Furthermore suppose that the value function $V(x)$ is continuous at $\bar x$ in direction $u$ if $u\not =0$ and the set-valued map $\Psi(x,y):=\Gamma-P(x,y)$ is  metrically subregular at $(\bar x, \tilde y,0)$ for all $\tilde y\in S(\bar x;u)$.  
If 
\begin{align*}
	\bigcup_{y\in  S(\bar x;u)}& 	\left \{\nabla_x P(\bar x, y)^T \lambda \left | \lambda \in \Lambda_u^0(\bar x,y; \mathcal C(\bar x,y;u)) \right .\right \}
	=\{0\},
\end{align*}
then $V(x)$ is Lipschitz around $\bar x$ in direction $u$, and
\begin{align*}
	\nonumber
	\emptyset &\not =\partial V(\bar x;u)
	\subseteq
	\bigcup_{y\in  S(\bar x;u)}  \left \{\nabla_x f(\bar x, y)+\nabla_x P(\bar x,  y)^T\lambda \left |\lambda \in \Lambda_u^1(\bar x,y; \mathcal C(\bar x,y;u))\right . \right \}.
\end{align*}
}
	\item[(iii)] Suppose that  there exists $\bar y\in S(\bar x)$ such that $S(x)$ is inner semicontinuous at $(\bar x,\bar y)$ in direction $u$  and the set-valued map $\Psi(x,y):=\Gamma-P(x,y)$ is  metrically subregular  at $(\bar x,\bar y,0)$. If 
	\begin{align*}
	&	 \left \{\nabla_x P(\bar x, \bar y)^T\lambda | \lambda \in \Lambda^0_u(\bar x,\bar y; \mathcal C(\bar x,\bar y;u))\cup \Lambda^0_0(\bar x,\bar y; \mathcal C(\bar x,\bar y;0)\cap\mathbb S)\right \}=\{0\},
	\end{align*}
	then $V(x)$ is Lipschitz around $\bar x$ in direction $u$, and 
	\begin{eqnarray*}
	\lefteqn{\emptyset \not =\partial V(\bar x;u) \subseteq}\\
	&&	 \left \{\nabla_x f(\bar x,\bar y)+\nabla_x P(\bar x, \bar y)^T\lambda|\lambda \in \Lambda^1_u(\bar x,\bar y; \mathcal C(\bar x,\bar y;u))\cup \Lambda^1_0(\bar x,\bar y; \mathcal C(\bar x,\bar y;0)\cap\mathbb S)\right \}.
 	\end{eqnarray*}
{
 	\item[(iv)] 
Suppose that  there exists $\bar y\in S(\bar x)$ such that $S(x)$ is inner calm at $(\bar x,\bar y)$ in direction $u$  and the set-valued map $\Psi(x,y):=\Gamma-P(x,y)$ is  metrically subregular  at $(\bar x,\bar y,0)$. 
If
 \begin{align*}
 	&	 \left \{\nabla_x P(\bar x, \bar y)^T\lambda | \lambda \in \Lambda^0_u(\bar x,\bar y; \mathcal C(\bar x,\bar y;u))\right \}=\{0\},
 \end{align*}
 then $V(x)$ is Lipschitz around $\bar x$ in direction $u$, and 
 \begin{eqnarray*}
 	\emptyset \not =\partial V(\bar x;u) \subseteq
 	 \{\nabla_x f(\bar x,\bar y)+\nabla_x P(\bar x, \bar y)^T\lambda|\lambda \in \Lambda^1_u(\bar x,\bar y; \mathcal C(\bar x,\bar y;u)) \}.
 \end{eqnarray*}
}
	\end{itemize}
\end{prop} 

For simplicity of notation we omit the equality constraints in the parametric optimization problem and consider the case where $\Gamma=\mathbb{R}^p_-$. 
Denote by $g(x,y):=P(x,y)$ and  $I_g(x,y):=\{i=1,\ldots,p|g_i(x,y)=0\}$. 
Since $\mathbb{R}^p_-$ is convex, by virtue of (\ref{convNormal}), we have  $$N_{\mathbb{R}_-^{p}}(g(x,y);\nabla g(x,y)(u,v))=N_{\mathbb{R}_-^{p}}(g(x,y))\cap \{\nabla g(x,y)(u,v)\}^\perp. $$
The critical cone is
\begin{eqnarray*}
\mathcal C(x,y;u)&=& \left \{v\in \mathbb{L}(x,y;u) |  V_-'(x;u)\leq   \nabla f(x,y)(u,v)\leq V_+'(x;u)\right \},
\end{eqnarray*} where 
$\mathbb{L}(x,y;u)=\{v|\nabla g_i(x,y)(u,v)\leq 0 \ \ (i\in I_g(x,y))\}.$

We define the set of the classical Lagrange multipliers and the singular Lagrange multipliers as
\begin{eqnarray*}
	\Sigma (x,y)& :=& \left \{\lambda\left |
	\nabla_y f(x,y) +\nabla_y g(x,y)^T\lambda=0,
	0\leq \lambda \perp g(x,y)\leq 0   \right. \right \},\\
	\Sigma^0(x,y)& :=& \left \{\lambda\left |
	\nabla_y g(x,y)^T\lambda=0,
	0\leq \lambda \perp g(x,y)\leq 0   \right. \right \},
\end{eqnarray*} respectively.
Now we can state  the following results based on Proposition \ref{smooth}. Unlike Bai and Ye \cite[Theorem 4.2]{BY} in which the RCR-regularity the Robinson stability are needed, the following results do not require the RCR-regularity and the Robinson stability.
\begin{prop} Let $u\in \mathbb{R}^n$. Consider the value function $V(x):=\inf_y \{ f(x,y)| g(x,y)\leq 0\}$ where $f$ and $g$ are smooth.
	\begin{itemize}
	\item[(i)]
	Suppose that  the restricted inf-compactness holds at $\bar x$ in direction $u$ with compact set $\Omega_u$. Furthermore suppose that $V(x)$ is continuous at $\bar x$ in direction $u$ if $u\not =0$ (e.g. when MFCQ holds at certain $\bar y\in S(\bar x)$, or the constraint mapping $g(x,y)$ is affine and the feasible region $\mathcal F(x)$ is nonempty near $\bar x$; {see e.g.,\cite[Proposition 4.1]{BY}}) and  that { the system $g(x,y)\leq 0$ is calm at $(\bar x, y)$ for all $y\in S(\bar x;u)$ (e.g. when MFCQ holds at  $(\bar x,y)$,} or the constraint mapping $g(x,y)$ is affine).  
	If 
	\begin{align*}
	\bigcup_{y\in  S(\bar x;u) \cap \Omega_u}
	&	\left(\left \{\nabla_x g(\bar x, y)^T\lambda_g\bigg|\begin{array}{l}
		\lambda_g \in \Sigma^0(\bar x,y)\cap \{\nabla g(\bar x, y)(u,v) \}^\perp,\\
		v\in \mathcal C(\bar  x, y;u) \end{array}\right \}\right. \\
	& \left.\bigcup \left \{\nabla_x g(\bar x, y)^T\lambda_g\bigg|\begin{array}{l}
		\lambda_g \in \Sigma^0(\bar x, y)\cap \{\nabla_y g(\bar x, y)v) \}^\perp,\\
		v\in \mathcal C({\bar x},y;0)\cap\mathbb S
	\end{array}
	\right \}\right)=\{0\},
	\end{align*}
	then $V(x)$ is Lipschitz around $\bar x$ in direction $u$, and
	\begin{align*}
	\nonumber
	\emptyset \not = &\partial V(\bar x;u) \subseteq \\
	&	\bigcup_{y\in  S(\bar x;u)\cap \Omega_u}\left ( \left \{\nabla_x f(\bar x, y)+\nabla_x g(\bar x,  y)^T\lambda_g \bigg|\begin{array}{l}
	\lambda_g \in \Sigma(\bar x,y)\cap \{\nabla g(\bar x, y)(u,v))\}^\perp,\\
	v\in \mathcal C(\bar  x,  y;u) \end{array}\right \} \right .\\
	& \left . \bigcup   \left \{\nabla_x f(\bar x, y)+\nabla_x g(\bar x,  y)^T\lambda_g\bigg|\begin{array}{l}
\lambda_g \in 	\Sigma(\bar x,y)\cap \{\nabla_y g(\bar x, y)v)\}^\perp ,\\
	v\in \mathcal C({\bar x},  y;0)\cap\mathbb S \end{array} 
	 \right \} \right ).
 	\end{align*}

	\item[(ii)]
Suppose that  $S(x)$ is inner calm* at $\bar x$ in direction $u$. {Furthermore suppose that $V(x)$ is continuous at $\bar x$ in direction $u$ if $u\not =0$} and the system $g(x,y)\leq 0$ is calm  at $(\bar x, y, 0)$ for all $y\in S(\bar x;u)$.  
If 
{\begin{align*}
	\bigcup_{y\in  S(\bar x;u)}
	&	\left(\left \{\nabla_x g(\bar x, y)^T\lambda_g\bigg|\begin{array}{l}
		\lambda_g \in \Sigma^0(\bar x,y)\cap \{\nabla g(\bar x, y)(u,v) \}^\perp,\\
		v\in \mathcal C(\bar  x, y;u) \end{array}\right \} \right)=\{0\},
\end{align*}
then $V(x)$ is Lipschitz around $\bar x$ in direction $u$, and
\begin{align*}
	\nonumber
	\emptyset \not = &\partial V(\bar x;u) \subseteq \\
	&	\bigcup_{y\in  S(\bar x;u)}\left ( \left \{\nabla_x f(\bar x, y)+\nabla_x g(\bar x,  y)^T\lambda_g \bigg|\begin{array}{l}
		\lambda_g \in \Sigma(\bar x,y)\cap \{\nabla g(\bar x, y)(u,v))\}^\perp,\\
		v\in \mathcal C(\bar  x,  y;u) \end{array}\right \}\right ).
\end{align*}}

	\item[(iii)] Suppose that  there exists $\bar y\in S(\bar x)$ such that $S(x)$ is inner semicontinuous at $(\bar x,\bar y)$ in direction $u$  and the  system $g(x,y)\leq 0$ is calm     at $(\bar x,\bar y)$. If 
	\begin{align*}
	&	\left \{\nabla_x g(\bar x, \bar y)^T\lambda_g\bigg|\begin{array}{l}
	\lambda_g \in \Sigma^0(\bar x,\bar y)\cap \{\nabla g(\bar x,\bar y)(u,v) \}^\perp,\\
	v\in \mathcal C(\bar  x, \bar y;u) \end{array}\right \} \\
	& \bigcup \left \{\nabla_x g(\bar x, \bar y)^T\lambda_g\bigg|\begin{array}{l}
	\lambda_g \in \Sigma^0(\bar x,\bar y)\cap \{\nabla_y g(\bar x,\bar y)v) \}^\perp,\\
	v\in \mathcal C({\bar x}, \bar y;0)\cap\mathbb S
	\end{array}
	\right \}=\{0\},
	\end{align*}
	then $V(x)$ is Lipschitz around $\bar x$ in direction $u$, and 
	\begin{align*}
	\nonumber
	\emptyset \not =\partial V(\bar x;u)\subseteq
	&	\left ( \left \{\nabla_x f(\bar x,\bar y)+\nabla_x g(\bar x, \bar y)^T\lambda_g \bigg|\begin{array}{l}
	\lambda_g \in \Sigma(\bar x,\bar y)\cap \{\nabla g(\bar x,\bar y)(u,v)\}^\perp ,\\
	v\in \mathcal C(\bar  x, \bar y;u) \end{array}\right \} \right .\\
	& \left . \bigcup   \left \{\nabla_x f(\bar x,\bar y)+\nabla_x g(\bar x, \bar y)^T\lambda_g\bigg|\begin{array}{l}
	\lambda_g \in \Sigma(\bar x,\bar y)\cap \{\nabla_y g(\bar x,\bar y)v\}^\perp ,\\
	v\in \mathcal C({\bar x}, \bar y;0)\cap\mathbb S \end{array} 
	 \right \} \right ).
 	\end{align*}
{
 	\item[(iv)] Suppose that  there exists $\bar y\in S(\bar x)$ such that $S(x)$ is inner calm at $(\bar x,\bar y)$ in direction $u$  and the  system $g(x,y)\leq 0$ is calm at $(\bar x,\bar y)$. If 
 \begin{align*}
 	&	\left \{\nabla_x g(\bar x, \bar y)^T\lambda_g\bigg|\begin{array}{l}
 		\lambda_g \in \Sigma^0(\bar x,\bar y)\cap \{\nabla g(\bar x,\bar y)(u,v) \}^\perp,\\
 		v\in \mathcal C(\bar  x, \bar y;u) \end{array}\right \} =\{0\},
 \end{align*}
 then $V(x)$ is Lipschitz around $\bar x$ in direction $u$, and 
 \begin{align*}
 	\nonumber
 	\emptyset \not =\partial V(\bar x;u)\subseteq
 	&	 \left \{\nabla_x f(\bar x,\bar y)+\nabla_x g(\bar x, \bar y)^T\lambda_g \bigg|\begin{array}{l}
 		\lambda_g \in \Sigma(\bar x,\bar y)\cap \{\nabla g(\bar x,\bar y)(u,v)\}^\perp ,\\
 		v\in \mathcal C(\bar  x, \bar y;u) \end{array}\right \}.
 \end{align*}
}
	\end{itemize}
\end{prop}

\subsection{Additive perturbations}
In this section we consider the case where the parametric optimization problem is obtained from additively perturbing an optimization problem,   i.e., $f(x,y):=f(y)$, $P(x,y):=x+P(y)$ where $f,P$ are  locally Lipschitz continuous. Consider the point $\bar x:=0$ and a directional perturbation in direction $u$, i.e., $x \xrightarrow{u} 0$  and we will give upper estimates  for the directional subdifferentials of the value function $V(x)$ at $\bar x=0$.

In this case, we define
\begin{eqnarray*}
\mathbb{L}(y;u){:=\mathbb{L}(\bar x,y;u)}&=&\{v|(u+ DP(y)(v))\cap T_\Gamma(P(y))\not =\emptyset\},\\
\mathcal C(y;u):=\mathcal C(\bar x,y;u)&=& \{ v\in \mathbb{L}(y;u) \left |f_-'(y;v)\leq V_+'(0;u),V_-'(0;u)  \leq f_+'(y;v) \right . \}.
\end{eqnarray*}

Let $y\in S(0)$ and $v\in \mathbb{R}^m$. {For $\alpha=0,1$, we define the set of generalized  Lagrange multipliers at $y$ in directions $\mathcal C(y;u)$ and $ \mathcal C(y;0)\cap\mathbb S$ as

\begin{eqnarray*}
\lefteqn{\Lambda_u^\alpha(y; \mathcal C(y;u)) }\\
&:=&  \left \{\lambda \left  |0\in \alpha\partial f(y;v)+ \partial \langle P,\lambda \rangle (y;v),\ \lambda \in 
 \bigcup_{d\in u+DP(y)(v)}  N_\Gamma(P(y); d ), \ \  v\in \mathcal C(y;u) \right . \right \},\\
\lefteqn{\Lambda_0^\alpha(y; \mathcal C(y;0)\cap\mathbb S)}\\
& :=& \left \{\lambda \left  |0\in \alpha\partial f(y;v)+ \partial \langle P,\lambda \rangle (y;v),\ \lambda \in  \bigcup_{d\in DP(y)(v)} N_\Gamma\left (P(y);d
\right ), \ v\in \mathcal C(y;0)\cap\mathbb S \right . \right \},
\end{eqnarray*}
  respectively. 
Then we have
\begin{eqnarray*}
M_u^\alpha(y;\mathcal C(y;u))& := &M_u^\alpha(\bar x,y;\mathcal C(y;u))=\left \{(\zeta,\lambda ) \left|\zeta=\lambda ,  \lambda\in \Lambda^\alpha_u(y; \mathcal C(y;u))\right. \right \},\\
M_0^\alpha(y;\mathcal C(y;0)\cap\mathbb S)& := & M_0^\alpha(\bar x,y;\mathcal C(y;0)\cap\mathbb S)=\left \{(\zeta,\lambda ) \left|\zeta=\lambda ,  \lambda\in \Lambda^\alpha_0 (y; \mathcal C(y;0)\cap\mathbb S)\right. \right \}.
\end{eqnarray*}}


First we show that the value function is directionally continuous at $\bar x=0$  provided that there is a solution $\bar y\in S(0)$ such that there is no  nonzero abnormal directional multipliers at $\bar y$. This extends the classical result of Gauvin and Dubeau \cite[Theorem 3.3]{GD} on the continuity of the value function to the directional continuity of the value function.
\begin{prop}\label{sufcont}
Consider the additive perturbed problem.  Suppose that the restricted inf-compactness condition holds at $\bar x=0$ in direction $u$. If there exists $\bar y\in S(0)$ satisfying 
\begin{align}
	\Lambda_u^0(\bar y; \mathcal C(\bar y;u))\cup \Lambda_0^0(\bar y; \mathcal C(\bar y;0)\cap\mathbb S)=\{0\},\label{singular}
\end{align}
then $V(x)$ is continuous around $\bar x=0$ in direction $u$.
\end{prop}
\beginproof
According to Proposition \ref{directional-LSC}, $V(x)$ is lower semicontinuous at $\bar x$ in direction $u$ under the directional restricted inf-compactness condition. Then it suffices to prove that $V(x)$ is upper semicontinuous at $\bar x$ in direction $u$, i.e., 
\[
\limsup_{x\xrightarrow{u}\bar x}V(x)\leq V(\bar x).
\] Since (\ref{singular}) holds,
there exists some direction $v$ such that
\begin{align*}
		0\in  \partial \langle P,\lambda \rangle (\bar y;v),\ \lambda\in N_\Gamma(P(\bar y);u+w) \mbox{ with } w\in DP(\bar y)(v), u+w\in T_\Gamma(P(\bar y))\implies \lambda=0
\end{align*}
 holds. Then by Proposition \ref{dMR}, the metric regularity for the set-valued map $\Psi(y):= \Gamma-P(y)$ holds at $(\bar y,0)$ in direction $(v, u)$ and so  there exist positive scalars $\delta,\varepsilon,\kappa$ such that for any  $(y,x)\in (\bar y,0)+\mathcal V_{\varepsilon,\delta}(v,u)$
\begin{equation*}
\mbox{dist}(y,\mathcal F(x))\leq\kappa\mbox{dist}(x+P(y),\Gamma).
\end{equation*}
Together with the lipschitzness of $P(y)$, this implies there exists $L>0$ such that
\begin{equation*}
	\mbox{dist}(\bar y,\mathcal F(x))\leq\kappa\mbox{dist}(x+P(\bar y),\Gamma)\leq\kappa\|x+P(\bar y)-P(\bar y)\|\leq\kappa L\|x-\bar x\|.
\end{equation*}
Then for any sequences $u^k\rightarrow u$ and $t_k\downarrow0$ satisfying that $$\lim_{k\rightarrow \infty} V(\bar x+t_ku^k)=\limsup_{x\xrightarrow{u}\bar x}V(x),$$ there exists $y^k\in\mathcal F(\bar x+t_ku^k)$ such that $\|y^k-\bar y\|\leq\kappa L\|x-\bar x\|$. Hence $y^k\rightarrow\bar y$. Consequently, $\limsup_{x\xrightarrow{u}\bar x}V(x)=\lim_{k\rightarrow \infty} V(\bar x+t_ku^k)\leq\lim_{k\rightarrow \infty} f(\bar x+t_ku^k,y^k)=f(\bar x,\bar y)=V(\bar x)$. This means $V(x)$ is upper semicontinuous at $\bar x=0$ in direction $u$. The proof is complete.
\endproof
Since the directional inner calmness* implies the directional local lower semicontinuity, from the proof it is easy to see that  Proposition \ref{sufcont} remains true if the directional restricted inf-compactness condition is replaced by the directional inner calmness*. Thus Theorem \ref{Prop3.4} applied to the additive perturbed problem has the following corollary.
\begin{cor}\label{cor3.1}  
	Consider the additive perturbed problem. Let $u\in \mathbb{R}^n$. 
	\begin{itemize}
	\item[(i)]
	Suppose that  the restricted inf-compactness holds at $\bar x=0$ in direction $u$ with compact set $\Omega_u$.  Moreover if $u\not =0$, suppose either $\Gamma$ {is geometrically derivable at $\tilde y$} for each $\tilde y\in S(0;u)\cap\Omega_u$ or $P(x,y)$ is directionally differentiable in direction $(u,v)$ for each $\tilde y\in S(0;u)\cap\Omega_u$ and $v\in\mathbb L(\tilde y;u)$.
	If 
	\begin{align*}
	& \bigcup_{y\in  S(0;u)\cap\Omega_u}  \Lambda_u^0(y; \mathcal C(y;u))\cup \Lambda_0^0(y; \mathcal C(y;0)\cap\mathbb S)=\{0\},
	\end{align*} 
	then $V(x)$ is Lipschitz around $\bar x=0$ in direction $u$ and
{	\begin{align*}
	\nonumber
	\emptyset \not =\partial V(0;u)\subseteq
	&  \bigcup_{y\in  S(0;u)\cap\Omega_u} \left( \Lambda_u^1(y; \mathcal C(y;u))\cup \Lambda_0^1(y; \mathcal C(y;0)\cap\mathbb S)\right).
 	\end{align*}}
{
\item[(ii)]
{If in {\rm (i)} the restricted inf-compactness at $\bar x$ in direction $u$ is replaced by the inner calmness* at $\bar x$ in direction $u$, then $\Lambda^\alpha_0(\bar x, y; \mathcal C(\bar x,y;0)\cap\mathbb S)(\alpha=0 ,1)$ can be removed in conditions of {\rm (i)} and one has,} 
if 
\begin{align*}
	& \bigcup_{y\in  S(0;u)}  \Lambda_u^0(y; \mathcal C(y;u))=\{0\},
\end{align*} 
then $V(x)$ is Lipschitz around $\bar x=0$ in direction $u$ and
\begin{align*}
	\nonumber
	\emptyset \not =\partial V(0;u)\subseteq
	&  \bigcup_{y\in  S(0;u)} \left( \Lambda_u^1(y; \mathcal C(y;u))\right).
\end{align*}
}
	\item[(iii)] Suppose that  there exists $\bar y\in S(\bar x)$ such that $S(x)$ is inner semicontinuous at $(0,\bar y)$. {Moreover if $u\not =0$, suppose either $\Gamma$ {is geometrically derivable at $\bar y$} or $P(y)$ is  directionally differentiable at $\bar y$ in direction $ v$ for each  $v\in\mathbb L(\bar y;u)$.} If 
	\begin{align*}
	&  \Lambda_u^0(\bar y; \mathcal C(\bar y;u))\cup \Lambda_0^0(\bar y; \mathcal C(\bar y;0)\cap\mathbb S)=\{0\},
	\end{align*}
	then $V(x)$ is Lipschitz around $\bar x=0$ in direction $u$ and 
	\begin{align*}
	\emptyset \not =\partial V(0;u)\subseteq
	&  \left( \Lambda_u^1(\bar y; \mathcal C(\bar y;u))\cup \Lambda_0^1(\bar y; \mathcal C(\bar y;0)\cap\mathbb S)\right).
	\end{align*}

	\item[(iv)] 
	{If in {\rm (iii)} the inner semicontinuity of $S(x)$ at $(\bar x,\bar y)$ in direction $u$ is replaced by the inner calmness at $(\bar x,\bar y)$ in direction $u$, then $\Lambda^\alpha_0(\bar x, \bar y; \mathcal C(\bar x,\bar y;0)\cap\mathbb S)(\alpha=0 ,1)$ can be removed in conditions of {\rm (i)} and one has,} 
	if  
\begin{align*}
	&  \Lambda_u^0(\bar y; \mathcal C(\bar y;u))=\{0\}
\end{align*}
then $V(x)$ is Lipschitz around $\bar x=0$ in direction $u$ and 
\begin{align*}
	\emptyset \not =\partial V(0;u)\subseteq
	&  \left( \Lambda_u^1(\bar y; \mathcal C(\bar y;u))\right).
\end{align*}

	\end{itemize}
\end{cor}
\beginproof (i)Since $P(x,y):=x+P(y)$, the Jacobian of $\nabla P(x,y)$ has full row rank, the NNAMCQ holds, hence the metric regularity of the set-valued map $\Psi(x,y):=\Gamma-x-P(y)$ holds automatically. If 
$$ \bigcup_{y\in  S(0;u)}  \Lambda_u^0(y; \mathcal C(y;u))\cup \Lambda_0^0(y; \mathcal C(y;0)\cap\mathbb S)=\{0\},$$ then 
by Proposition \ref{sufcont},
 the value function is continuous in direction $u$. {Taking into account Proposition \ref{Prop3.4}, (ii)-(iv) can be proved following a similar process.}
\endproof

Consider the following optimization problem
\begin{eqnarray*}
 (\widetilde{P})~~~~\quad \min_{y\in \Omega} &&f(y) \quad \mbox{ s.t. }  g(y)\leq 0, h(y)=0 
\end{eqnarray*}
where $f:\mathbb R^{m}\rightarrow\mathbb R,\ g:\mathbb R^{m}\rightarrow\mathbb R^p, h:\mathbb R^{m}\rightarrow\mathbb R^q$ are locally Lipschitz continuous and $\Omega\subseteq\mathbb R^p$ is a closed set,
{and its additive perturbation
\begin{eqnarray*}
(\widetilde{P})_{s,t}~~~~\quad \min_{y\in \Omega} &&f(y) \quad \mbox{ s.t. } g(y)+s \leq 0, h(y)+t=0,
\end{eqnarray*}
where $s\in \mathbb{R}^p, t\in \mathbb{R}^q$ are the parameters. The value function is
$$V(s,t):=\inf_{y\in \Omega} \{ f(y) | g(y)+s\leq 0, h(y)+t=0\}.$$
Given a direction $u:=(\eta,\beta) \in \mathbb{R}^{p+q}$, we consider a directional perturbation in direction $u$, i.e., $x:=(s,t) \xrightarrow{u} (0,0)$. 
In this case, 
\begin{eqnarray*}
	&& x+P(y)=(s+g(y),t+h(y),y),\ \Gamma=\mathbb R^p_-\times\{0\}^q\times \Omega,\\
	&& T_\Gamma(P(y))= T_{\mathbb R_-^p}(g(y))\times T_{\{0\}^q}(h(y)) \times T_\Omega(y),\\ 
	&& \mathbb{L}(y;u)=\{v|\{\eta_i+Dg_i(y)(v)\}\cap\mathbb R_-\neq\emptyset(i\in I_g(y)),0\in\beta+Dh(y)(v), v\in T_\Omega(y)\},\\
	&& \mathcal C(y;u)= \{ v\in \mathbb{L}(y;u) \left |f_-'(y;v)\leq V_+'(0;u),V_-'(0;u)  \leq f_+'(y;v) \right . \}.
	\end{eqnarray*}
	Moreover by \cite[Proposition 3.2]{YZ17}, for any $(\eta,\beta,0)+(\xi,\zeta,v))\in T_\Gamma(P(y))$, we have
$$ {N_\Gamma(P(y);(\eta,\beta,0)+(\xi,\zeta,v))}=
N_{{\mathbb R}^p_-}(g(y);\eta+\xi)\times N_{\{0\}^q}(h(y);\beta+\zeta)\times N_\Omega(y;v),
$$ 
 Note that since $T_{\{0\}^q}(h(y))=\{0\}^q$, we have $N_{\{0\}^q}(h(y);\beta+\zeta)=\{0\}^q$ and  by (\ref{convNormal}), we have
$N_{\mathbb R_-^p}(g(y);\eta+\xi)=N_{\mathbb R_-^p}(g(y))\cap\{\eta+\xi\}^\perp$ for $\eta+\xi\in T_{\mathbb R_-^p}(g(y)).$ Let $y\in S(0)$ and $v\in \mathbb{R}^m$. For $\alpha=0,1$,  we define the set of generalized Lagrange multipliers at $y$  in directions $\mathcal C(y;u)$ and $\mathcal C(y;0)\cap\mathbb S$ as
\begin{eqnarray*}
	\lefteqn{\Lambda^\alpha_u(y; \mathcal C(y;u))}\\
	& :=& \left \{(\lambda_g,\lambda_h)\left  |\begin{array}{ll}
		0\in \alpha\partial f(y;v)+ \partial \langle g ,\lambda_g\rangle (y;v)+\partial \langle h ,\lambda_h\rangle (y;v)+\displaystyle  N_\Omega\left (y;v\right ),v\in \mathcal C(y;u),  \\
		0\leq\lambda_g\perp g(y), {\lambda_g \in \{\eta+\xi\}^\perp, \ \mbox{for some}\ \xi\in Dg(y)(v)\ \mbox{with}\ \eta+\xi\in T_{\mathbb R^p_-}(g(y))}
	\end{array}\right .\right \},\\
	\lefteqn{\Lambda^\alpha_0(y; \mathcal C(y;0)\cap\mathbb S)}\\& :=& \left \{(\lambda_g,\lambda_h)\left  |\begin{array}{ll}
		0\in \alpha\partial f(y;v)+ \partial \langle g ,\lambda_g\rangle (y;v)+\partial \langle h ,\lambda_h\rangle (y;v)+   N_\Omega\left (y;v\right ),v\in \mathcal C(y;0)\cap\mathbb S \\
		0\leq\lambda_g\perp g(y), {\lambda_g \in \{\xi\}^\perp, \ \mbox{for some}\ \xi\in Dg(y)(v)\cap T_{\mathbb R^p_-}(g(y))}
	\end{array}\right.\right\},
\end{eqnarray*} respectively. 
When $u=0$, $\mathcal C(y;0)\cap\mathbb S\subseteq \mathcal C(y;0)$. Since $0\in \mathcal C(y;0)$, {taking $u=0$ in the above},  the directional multipliers become the nondirectional multipliers and we have
\begin{eqnarray*}
\lefteqn{\Lambda_0^1(y; \mathcal C(y;0))\cup \Lambda_0^1(y; \mathcal C(y;0)\cap \mathbb{S})}\\
&&= \left \{(\lambda_g,\lambda_h)|0\in  \partial f(y)+ \partial\langle g ,\lambda_g\rangle (y)+\partial \langle h,\lambda_h\rangle (y)+N_\Omega\left (y\right ), 0\leq\lambda_g\perp g(y)   \right \},\\
\lefteqn{\Lambda_0^0(y; \mathcal C(y;0))\cup \Lambda_0^0(y; \mathcal C(y;0)\cap\mathbb S)}\\
&&= \left \{(\lambda_g,\lambda_h)|0\in  \partial\langle g ,\lambda_g\rangle (y)+\partial \langle h,\lambda_h\rangle (y)+N_\Omega\left (y\right ), 0\leq\lambda_g\perp g(y)   \right \}.
\end{eqnarray*}
In this case taking the convex hull, Corollary \ref{cor3.1}(i) recovers the classical result in Clarke \cite[Corollary 1 of Theorem 6.52]{Clarke}.}

\subsection{Danskin's Theorem}

	Recall that the classical Danskin's theorem 
	can be stated as follows.
\begin{thm}[Danskin's Theorem]\cite[Problem 9.13, Page 99]{Danskin}
  Let $f:\mathbb{R}^{n+m}\rightarrow \mathbb{R}$ be continuous and $\Gamma$ be a compact subset of $\mathbb{R}^m$. Suppose that for a given neighborhood $U(\bar x)$ of $\bar x$, the gradient $\nabla_x f(x,y)$ exists and is continuous (jointly) for $(x,y)\in U(\bar x)\times \Gamma$. Then the value function $V(x):=\inf_{y\in \Gamma} f(x,y)$ is Lipschitz around $\bar x$ and
	\begin{align*}
	\partial^c V(\bar x)={\rm  co } \left \{\nabla_x f(\bar x, y) | {y\in  S(\bar x)}  \right \}.
	\end{align*}
\end{thm}
In the following, we introduce the directional Danskin's theorem, where the compactness of $\Gamma$ is replaced by the weaker condition, directional restricted inf-compactness condition, and sharper estimates for both directional limiting and Clarke subdifferentials of value function are obtained.
\begin{prop}[Directional Danskin's Theorem]  Let $f:\mathbb{R}^{n+m}\rightarrow \mathbb{R}$ be continuous and $\Gamma$ is closed. Suppose that for a given neighborhood $U(\bar x)$ of $\bar x$, the gradient $\nabla_x f(x,y)$ exists and is continuous (jointly) for $(x,y)\in U(\bar x)\times \Gamma$. 
	\begin{itemize}
		\item[{\rm (i)}]
		Suppose the restricted inf-compactness holds at $\bar x=0$ in direction $u$ with compact set $\Omega_u$.   
		Then $V(x):=\inf_{y\in \Gamma} f(x,y)$ is Lipschitz around $\bar x$ in direction $u$ and
		\begin{align}
		& 	\emptyset \not =\partial V(\bar x;u)\subseteq
				\bigcup_{y\in  S(\bar x;u)\cap\Omega_u}  \left \{\nabla_x f(\bar x, y)  \right \},\label{inclusion}\\
	&		\partial^c V(\bar x;u)={\rm  co } \left \{\nabla_x f(\bar x, y) | {y\in  S(\bar x;u)\cap\Omega_u}  \right \}.\label{equality}
		\end{align}
		\item[{\rm (ii)}] Suppose that there exists $\bar y\in S(\bar x)$ such that $S(x)$ is inner semicontinuous at $(\bar x,\bar y)$ in direction $u$. Then $V(x)$ is Lipschitz around $\bar x$ in direction $u$ and 
		\begin{align*}
			\nonumber
			\partial V(\bar x;u)=
			& \left \{\nabla_x f(\bar x,\bar y)\right \}.
		\end{align*}
	\end{itemize}
\end{prop}
\beginproof  
(i) According to Proposition \ref{directional-LSC}, $V(x)$ is lower semicontinuous at $\bar x$ in direction $u$ under the directional restricted inf-compactness condition. Take any sequences $u^k\rightarrow u$ and $t_k\downarrow0$ satisfying that $\lim_{k\rightarrow \infty} V(\bar x+t_ku^k)=\limsup_{x\xrightarrow{u}\bar x}V( x)$. We have for any $\bar y\in S(\bar x;u)$,
$$\limsup_{x\xrightarrow{u}\bar x}V( x)=\lim_{k\rightarrow \infty} V(\bar x+t_ku^k)\leq\lim_{k\rightarrow \infty} f(\bar x+t_ku^k,\bar y)=f(\bar x,\bar y)=V(\bar x),$$ which means that $V(x)$ is upper semicontinuous at $\bar x$ in direction $u$. Hence $V(x)$ is continuous at $\bar x$ in direction $u$.

{Next we prove  (\ref{inclusion}).
(i) 
Let $\zeta\in\partial V(\bar x;u)$. Then by Definition \ref{ads}, there exist sequences $t_k\downarrow0,\ u^k\rightarrow u,\ \zeta^k\rightarrow\zeta$ such that $V(\bar x+t_ku^k)\rightarrow V(\bar x)$
and 
$\zeta^k\in \widehat \partial V(\bar x+t_ku^k)$. It follows that 
$V(\bar x+t_ku^k)<V(\bar x)+\varepsilon$ for all $k$ large enough and hence by  the directional restricted inf-compactness, there exists $y^k\in S(\bar x+t_ku^k)\cap \Omega_u$.  Passing to a subsequence if necessary, we may assume that $y^k\rightarrow\tilde y$. Hence $\tilde y\in S(\bar x;u)\cap \Omega_u$.\\
For each $k$, since $\zeta^k\in \widehat \partial V(\bar x+t_ku^k)$, there exists a neighborhood ${\mathcal U}^k$ of $\bar x+t_ku^k$ satisfying
\begin{equation*}
	V(x)-V(\bar x+t_ku^k)-\langle\zeta^k,x-(\bar x+t_ku^k)\rangle+\frac{1}{k}\|x-(\bar x+t_ku^k)\|\geq0\ \forall x\in{\mathcal U}^k.
\end{equation*}
It follows from the fact $V(x)=\displaystyle \inf_{y} \left \{f(x,y)+\delta_{ C}(y) \right \}$ and $y^k\in S(\bar x+t_ku^k)$, that
\[
f(x,y^k)-\langle\zeta^k,x-(\bar x+t_ku^k)\rangle+\frac{1}{k}\|x-(\bar x+t_ku^k)\|\geq f(\bar x+t_ku^k,y^k),
\]
for any $x\in{\mathcal U}^k$.
Hence, the function 
$$
x\rightarrow f(x,y^k)-\langle\zeta^k,x-(\bar x+t_ku^k)\rangle+\frac{1}{k}\|x-(\bar x+t_ku^k)\|$$ attains its local minimum at $x=\bar x+t_ku^k$. Then by the well known Fermat's rule (see e.g., \cite[Proposition 1.30(i)]{Mor}) and the calculus rule in Proposition \ref{Calculus}, 
\begin{equation}
	0\in \nabla_x f(\bar x+t_ku^k,y^k)-\zeta^k+\frac{1}{k} {\mathcal B}.
	\label{boundedn}
\end{equation}
 Since $\nabla_xf(x,y)$ is continuous (jointly) on $ U(\bar x)\times C$, taking the limit of (\ref{boundedn}) as $k\rightarrow\infty$, we obtain
\begin{equation*}
	0\in \nabla_x f(\bar x,\tilde y)-\zeta.
\end{equation*} Inclusion (\ref{inclusion}) is proved.}

We now prove  equality (\ref{equality}). Consider an arbitrary element $\tilde y\in S(\bar x;u)\cap\Omega_u$. Then there exist sequences $t_k\downarrow0, u^k\rightarrow u$ and $y^k\rightarrow\tilde y$ with $y^k\in S(\bar x+t_ku^k)\cap\Omega_u$. By Taylor expansion, one has 
\begin{align*}
	-V(x)&\geq -f(x, y^k)\\
	&=-f(\bar x+t_ku^k,y^k)-\nabla_xf(\bar x+t_ku^k,y^k)^T(x-(\bar x+t_ku^k))+o(\|x-(\bar x+t_ku^k)\|)\\
	&=-V(\bar x+t_ku^k)-\nabla_xf(\bar x+t_ku^k,y^k)^T(x-(\bar x+t_ku^k))+o(\|x-(\bar x+t_ku^k)\|).
\end{align*} Hence $-\nabla_xf(\bar x+t_ku^k,y^k)\in\widehat \partial (-V)(\bar x+t_ku^k)$. Since $-V(x)$ is {Clarke regular (\cite[Page 99]{Clarke})}, we have $-\nabla_xf(\bar x+t_ku^k,y^k)\in\partial^c(-V)(\bar x+t_ku^k)$. By the continuity of $\nabla_xf(x,y)$ and (\ref{clarke}), taking the limit as $k\rightarrow\infty$, one has $-\nabla_xf(\bar x,\tilde y)\in \partial^c(-V)(\bar x;u)$. By the choice of $\tilde y$, one has
\begin{align*}
	\partial^c (-V)(\bar x;u)\supseteq
	\bigcup_{y\in  S(\bar x;u)\cap\Omega_u}  \left \{-\nabla_x f(\bar x, y)  \right \}.
\end{align*}
On the other hand, (\ref{inclusion}) implies that
\begin{align*}
	& 	\partial^c (-V)(\bar x;u)\subseteq
	{\rm co}\left\{\bigcup_{y\in  S(\bar x;u)\cap\Omega_u}  \left \{-\nabla_x f(\bar x, y)  \right \}\right\}.
\end{align*}
In summary, one obtains
\begin{align*}
	& 	\partial^c (-V)(\bar x;u)=
	{\rm co}\left\{\bigcup_{y\in  S(\bar x;u)\cap\Omega_u}  \left \{-\nabla_x f(\bar x, y)  \right \}\right\}.
\end{align*}
By (\ref{clarke}), $\partial^c (-V)(\bar x;u)=-\partial^cV(\bar x;u)$. Hence,
\begin{align*}
	& 	\partial^c V(\bar x;u)=
	{\rm co}\left\{\bigcup_{y\in  S(\bar x;u)\cap\Omega_u}  \left \{\nabla_x f(\bar x, y)  \right \}\right\}.
\end{align*}
(ii) Let $\tilde y=\bar y$ in the above proof. One can obtain
\begin{align*}
	& 	\partial^c V(\bar x;u)=
	 \left \{\nabla_x f(\bar x, \bar y)  \right \}.
\end{align*}
\endproof

\vskip 6mm
\noindent{\bf Acknowledgements}

\noindent
The authors would like to thank the anonymous referees for their helpful suggestions and comments. The research of the first author was partially supported by Hong Kong Research Grants Council PolyU153036/22p. The research of the second author was partially supported by NSERC.


\end{document}